\documentclass[mnsc,nonblindrev]{informs3} % current default for manuscript submission

\OneAndAHalfSpacedXI % current default line spacing
\usepackage{natbib}
 \bibpunct[, ]{(}{)}{,}{a}{}{,}%
 \def\bibfont{\small}%

	\usepackage{amsmath}
	\usepackage{graphicx}
	\usepackage{enumerate}
	\usepackage{xcolor}
	\usepackage{natbib} %comment out if you do not have the package
	\usepackage{url} % not crucial - just used below for the URL
	\usepackage{multirow}
	\usepackage{color}
        \usepackage{soul}
        \usepackage{parallel}
        \usepackage{multicol}
        \usepackage{rotfloat} 
        \usepackage{booktabs}
        \usepackage{caption}
        \usepackage{subcaption}
        \usepackage{lscape}
	\usepackage{amssymb}
	\usepackage{multirow}
 \usepackage{wrapfig,lipsum,booktabs}

 \usepackage{algorithm}
\usepackage{algpseudocode}

\TheoremsNumberedThrough     % Preferred (Theorem 1, Lemma 1, Theorem 2)
\ECRepeatTheorems

\EquationsNumberedThrough    % Default: (1), (2), ...
\MANUSCRIPTNO{}

\renewcommand{\theARTICLETOP}{}

\begin{document}
%%%%%%%%%%%%%%%%

%\setcounter{page}{1}
%\YEAR{0000}% e.g.,, 2005
%\DOI{10.1287/xxxx.0000.0000}%


\RUNTITLE{Routing Strategies for Equitable Housing Coordination}

% \TITLE{Bundling Information Goods of Decreasing Value}
\TITLE{Queue Routing Strategies to Improve Equitable Housing Coordination in New York City}


% commented out the author names for blind review

\ARTICLEAUTHORS{%
\AUTHOR{Yaren Bilge Kaya }
\AFF{Industrial Engineering Department, Northeastern University, Boston, USA  \EMAIL{kaya.y@northeastern.edu}} %, \URL{}}
\AUTHOR{and Kayse Lee Maass}
\AFF{Industrial Engineering Department, Northeastern University, Boston, USA  \EMAIL{k.maass@northeastern.edu}}
} % end of the block

\ABSTRACT{%
Runaway and homeless youth (RHY) are a group of youth and young adults who are at high risk of being exploited through human trafficking. Although access to housing and support services is an effective way to decrease their vulnerability to being exploited, research reveals that coordination of these services provided to RHY by non-profit and government organizations is neither standardized, nor efficient. This situation often causes decreased, delayed, and inequitable access to these scarce housing resources. In this study, we aim to increase the housing system efficiency and reduce the barriers that are contributing to inequitable access to housing through simulation modeling and analyses. Specifically, we simulate a set of crisis and emergency shelters in New York City, funded by a single governmental organization, considering a queuing network with pools of multiple parallel servers, servers with demographic eligibility criteria, stochastic RHY arrival, impatient youth behaviour (possibility of abandonment), and a decision-maker (coordinator) that determines which server pool RHY is routed to. This simulation allows us to evaluate the impact of different queue routing strategies. Our simulation results show that by changing the way RHY is routed to shelters, we can reduce the average wait time by approximately a day and decrease the proportion of RHY abandoning the shelters by 13\%.

}%

%\KEYWORDS{deterministic inventory theory; infinite linear programming duality;
%  existence of optimal policies; semi-Markov decision process; cyclic schedule}


\maketitle
%%%%%%%%%%%%%%%%%%%%%%%%%%%%%%%%%%%%%%%%%%%%%%%%%%%%%%%%%%%%%%%%%%%%%%

% Samples of sectioning (and labeling) in GNNSFC
% NOTE: (1) \section and \subsection do NOT end with a period
%       (2) \subsubsection and lower need end punctuation
%       (3) capitalization is as shown (title style).
%
%\section{Introduction.}\label{intro} %%1.
%\subsection{Duality and the Classical EOQ Problem.}\label{class-EOQ} %% 1.1.
%\subsection{Outline.}\label{outline1} %% 1.2.
%\subsubsection{Cyclic Schedules for the General Deterministic SMDP.}
%  \label{cyclic-schedules} %% 1.2.1
%\section{Problem Description.}\label{problemdescription} %% 2.

% Text of your paper here

%To get the document ready for arxiv submission uncomment these:
\setlength{\footskip}{24pt}
\pagestyle{plain}

\section{Introduction}

Homelessness is a growing humanitarian problem around the world \citep{Henry-2021}, with an estimated 4.2 million youth and young adults in the United States (US) experiencing homelessness annually \citep{ChapinHall-2018}. Research indicates that runaway and homeless youth (RHY) are more susceptible to physical and emotional trauma, sexual abuse, drug misuse, mental health problems, and illnesses before and during their homelessness period when compared to their housed counterparts \citep{Martijn-2006, Heerde-2015, Mallett-2005}. These factors, coupled with a lack of community support, are likely to render them more vulnerable to precarious situations and human trafficking \citep{Greenbaum-2017}. 

Exiting homelessness is challenging \citep{Morton-2018}. RHY encounter systemic, situational, and interpersonal barriers across multiple systems while trying to exit homelessness \citep{Sample-2020}. However, findings indicate that stable housing, economic resources, and physical and mental health support is crucial to break the cycle of homelessness \citep{Zlotnick-1999, Thompson-2004}. In the US, these housing and support services are often provided by governmental and non-profit organizations through limited funds \citep{Nelson-1994}. Consequently, in most US communities, the demand for these resources significantly surpasses the supply, and lack of coordination contributes to inefficiencies and inequitable access \citep{Clawson-2009}.  

The US Department of Housing and Urban Development launched the Coordinated Entry System (CES) to address the inefficiencies within the public housing system and CES was nationally-recognised shortly after its launch \citep{HUD-2015, Crossroads-2023}. With CES the coordinators aimed for increased access, reduced barriers for unhoused persons (clients), strategically prioritized resources, and consistency across the system. To achieve these objectives, CES replaced individually managed wait-lists and first come first serve (FCFS) application processes with a person-centered, standardized tool. The new system enhanced access and referrals by adopting a \textit{``No-Wrong Door Access"} approach and ensuring shared tools and processes across organizations. 

Despite the CES's significance as an initial effort to enhance housing provision, studies show that individually done client assessments and client ``hand-offs`` between service providers may be contributing to inequitable housing access \citep{Balagot-2019, Ecker-2022}. Moreover, studies reveal that CES necessitates improvements related to technology and user education, as well as federal policy modifications \citep{Hornung-2020}. Presently, there is little evidence that coordinated systems improve individual-level outcomes such as length of stay in housing \citep{Dickson-gomez-2020, Ecker-2022}. 

To address these observed gaps, in this trans-disciplinary study, we seek to improve equitable access to scarce public housing resources through simulation modeling and analyses. We simulate a set of RHY crisis-emergency shelters that are funded by a single governmental agency, taking into account factors such as (i) exponential inter-arrival times, (ii) generally distributed service times, (iii) multiple server pools with demographic eligibility criteria, (iv) impatient clients with generally distributed patience times, (v) heterogeneous server and customers, and (vi) a decision-maker (coordinator) that determines which server pool RHY is routed to. Our study population is RHY between 16-24 years of age in New York City (NYC), where the rate of homelessness is the highest in the US \citep{Morton-2019}. The simulation enables us to introduce and compare various queue routing strategies and evaluate outcomes of these strategies based on different demographic characteristics such as age, gender, sexual orientation, immigration status, and human trafficking victimization status. To the best of our knowledge this study is the first attempt in the operations research and analytic literature to compare different interpretible routing strategies for improving the housing allocation process in terms of efficiency and equity.

The remainder of the paper is structured as follows. We provide a review of the literature in Section 2; present the routing strategies and service provision at shelters in Section 3. In Section 4, we present our computational setup; in Section 5, our  ensuing results, and insights regarding service quality, and equitable access. We conclude by summarizing our contributions, limitations and future research directions in Section 6.

\section{Literature Review}
The application of operations research (OR) and analytics techniques can provide solutions to various challenges in public housing systems.  Previous research related to OR and analytics efforts to reduce homelessness and improve access to housing resources largely focus on three broad areas: (i) location planning \citep{Johnson-2000, Johnson-2007, Maass-2020}, (ii) capacity planning \citep{Kaya-2022a, Kaya-2022b, Miller-2022}, and (iii) allocation mechanisms or routing strategies \citep{Kaplan-1986, Chan-2018, Azizi-2018, Arnosti-2020, Rahmattalabi-2022, Kaya-2023}. 

In this study, we model the allocation of crisis-emergency homeless shelter beds as a queuing problem and compare different routing strategies. Our work draws upon and contributes to two primary areas: (i) OR for improving equitable access to scarce public resources and (ii) queue routing strategies.

\subsection{OR and Housing Allocation Mechanisms}
In this paper, we study the problem of designing equitable and interpretable policies that effectively match heterogeneous RHY to scarce housing resources. To consider different types of housing resources, \cite{Rahmattalabi-2022} model this problem as a multi-class multi-server queuing system, where each individual is assigned to a queue where they wait to be matched to a resource considering the eligibility structure. The authors present a methodology based on causal inference and propose an eligibility structure that accounts for Non-stress test vulnerability score and demographic information such as age and race \citep{Rahmattalabi-2022}. In an earlier study, \cite{Kaplan-1986} formulates the allocation of affordable housing as a queuing problem and studies waiting times and development diversity under various priority rules. Although both these studies provide innovative approaches, the authors recognize the need for further exploration and evaluation of alternative policies. 

\cite{Chan-2018} identifies inefficiencies in the current housing system by pointing out the assignment decisions that are made manually by humans working in the housing communities. To make strides towards an autonomous intelligent agent instead of more manual work, \cite{Chan-2018} formulate the problem of assigning homeless youth to housing programs subject to resource constraints as a multiple multi-dimensional knapsack problem. The authors propose the \textit{``Minimum Resources Consumed First"}  (MRCF) strategy that assigns youth that consume the least amount of resources first. While this study provides an efficient mechanism for housing allocation, it neglects considerations of fairness and equity from the perspective of homeless youth.

\cite{Azizi-2018} examine the allocation of permanent supportive housing and rapid rehousing, emphasizing the importance of fairness and equity measures. The authors propose a data-driven mixed-integer optimization formulation that is flexible and takes into account demographic characteristics such as age and history of substance use. Their strategy, which we define as the ``Most Likely to Exit Homelessness First" (MLEHF), focuses on maximizing the probability of unhoused youth's safe and stable exit from the housing program. \cite{Kaya-2023} take a different approach and study equitable access for the most vulnerable RHY, considering demographic characteristics including human trafficking history, mental health or substance abuse issues, and immigrant status. The authors present the \textit{``Most Vulnerable First"} (MVF) strategy that prioritizes youth who are at greatest risk of experiencing human trafficking.  They propose expanding the current housing capacity and utilizing priority thresholds to improve equitable access. 

Although the prior literature on housing allocation mechanisms draws upon and incorporates queue routing strategies, the broader queuing literature contains additional routing strategies that can be applied to housing access and CES. Therefore, our study presents a comparative analysis of various allocation mechanisms to address this gap in the housing assignment literature.

\subsection{Queue Routing Strategies}
Queuing theory literature comprises numerous empirical studies that explore various routing strategies. While some of these studies investigate quality-of-service, fairness, and equity metrics from the customers' standpoint, others adopt a server-centric approach. The selection of these perspectives is contingent upon the service environment and the characteristics of the customers. In this section, we present a range of queue routing strategies and concentrate specifically on those that can be deemed appropriate as public housing allocation mechanisms.

%Queuing theory literature contains several empirical studies on routing strategies. Some of these studies consider quality-of-service, fairness, and equity measures from the customers' perspective and some choose to investigate them through the servers' point of view. This choice highly depends on the service setting and customer characteristics. In this section, we introduce various routing strategies and particularly focus on those that are suitable to be considered as public housing allocation mechanisms.

Numerous studies have explored queuing models where customers choose between multiple parallel servers, with some focusing on the \textit{``Shortest-Queue First"} (SQF) approach \citep{Lehtonen-1984, Whitt-1986, Sparaggis-1995, Mukherjee-2016}. This common strategy assigns incoming customers to the server or server pool with the shortest queue \citep{Ephremides-1980, Johri-1989, Hordijk-1990, Koole-1999, Akgun-2011}. In an early study, \cite{Ephremides-1980} introduces SQF and consider a queuing model in which arriving customers have to choose between $m$ parallel \textit{homogeneous} servers. Later studies expand this approach for different arrival and server behaviours. For example, \cite{Johri-1989} incorporates exponential servers with state-dependent service rates. \cite{Akgun-2011} expand the SQF for arbitrary arrivals where only a subset are flexible, multiple-server stations, and abandonment. Consequently, in this study, we also evaluate the effectiveness of the SQF strategy for the NYC crisis and emergency housing system.

%Many studies consider queuing models in which arriving customers have to choose between m parallel servers \citep{Lehtonen-1984, Whitt-1986, Sparaggis-1995, Mukherjee-2016}. A common approach in these type of systems is known to be the \textit{``Shortest-Queue First"} (SQF), where customers are assigned to the server or server pool with the shortest queue \citep{Ephremides-1980, Johri-1989, Hordijk-1990, Koole-1999, Akgun-2011}. In an early study, \cite{Ephremides-1980} introduces SQF and consider a queuing model in which arriving customers have to choose between m parallel \textit{homogeneous} servers. Later studies expand this approach for different arrival and server behaviours. For example, \cite{Johri-1989} incorporates exponential servers with state-dependent service rates. \cite{Akgun-2011} expand the SQF for flexible arrivals where arbitrary arrivals where only a subset are flexible, multiple-server stations, and abandonment. Consequently, in this study, we also evaluate the effectiveness of the SQF strategy for NYC crisis and emergency housing system.

\cite{Armony-2005} introduces the \textit{``Fastest Server Available First"} (FSF) policy which aims to minimize the steady-state queue length and virtual waiting time in call centers. Although it provides useful insights, this policy fails to consider fairness from the server's perspective by putting additional burden on the fastest servers. Therefore, it is not recommended to be used in public service provision settings. Later, \cite{Armony-2011} expand this study by incorporating impatient customer behaviour where each customer has a maximum amount of time they are willing to wait to receive service. They illustrate the importance of incorporating queue abandonment accurately while modeling service systems. Motivated by this, in our study, we consider RHY as individuals with limited patience. 

In an attempt to achieve fairness toward servers, \cite{Atar-2011} and \cite{Adan-2012} recommend the \textit{``Longest-Idle Server Pool First"}  (LISF) policy, which is commonly used in call centers and deemed fair. While doing so, \cite{Atar-2011} model the system as a single queue with a fixed number of server pools of heterogeneous exponential servers and routes a customer to the pool with the longest cumulative idleness among pools. \cite{Adan-2012} expand this model by incorporating skill based service (where servers can serve a subset of customer types). Drawing upon these aforementioned studies, we model the NYC crisis and emergency housing system as a queuing network with multiple heterogeneous server pools, which represent different shelters with inclusion and exclusion criteria, and homogeneous servers that represent beds within a shelter (i.e., homogeneous servers within heterogeneous server pools). %Drawing upon these aforementioned studies, we model the NYC crisis-emergency housing system considering a queuing network with multiple heterogeneous server pools (to reflect the different inclusion and exclusion criteria among shelters) of homogeneous servers (beds within a shelter). 
Furthermore, in our study, we are interested in looking at fairness from both the customer's and server's perspectives. %Furthermore, we incorporate skill based service approach since shelters have inclusion and exclusion criteria that dictate who may receive services at the shelter. 

Focusing on a healthcare setting, \cite{Mandelbaum-2012} model the transition between an emergency department and internal wards as a single centralized queue and multiple server pools (a.k.a. the inverted-V model). Resulting in the same server fairness as LISF, the authors introduce the \textit{``Randomized Most Idle"} (RMI) strategy, which assigns a customer to an available pool with probability equal to the fraction of idle servers in that pool out of the total number of idle servers in the system \citep{Mandelbaum-2012}. Compared to LISF, the RMI policy requires less real-time information from the system, which makes it easier to implement. 

In addition to these, there are various other queue routing strategies that are introduced in the analytics literature such as \textit{``Shortest Idle Server First"} (SISF) \citep{Madadi-2023} and \textit{``Minimum Expected Delay Faster Server First"} \citep{Tezcan-2008}. However, we do not include those in this study due to their irrelevance to our aim. In Table \ref{t:literature}, we present the summary of routing strategies explained in our literature review with their abbreviations and corresponding key references. 

\begin{table}[h]
\centering
\caption{The relevant CES queue routing strategies introduced in the analytics literature so far, presented with associated abbreviations and references.}
\begin{tabular}{lcc}
\toprule
\textbf{Routing   Strategy} & \multicolumn{1}{l}{\textbf{Abbreviation}} & \textbf{Study} \\ 
\midrule
\multirow{2}{*}{Fastest Server Available First} & \multirow{2}{*}{FSF} & \cite{Armony-2005} \\
& & \cite{Armony-2011}\\
 \midrule
\multirow{2}{*}{Longest-Idle Server Pool First} & \multirow{2}{*}{LISF} & \cite{Atar-2011} \\
 &  & \cite{Adan-2012} \\
 \midrule
Minimum Resources Consumed First & MRCF & \cite{Chan-2018} \\
\midrule
\multirow{2}{*}{Most Likely to Exit Homelessness First} & \multirow{2}{*}{MLEHF} & \cite{Rahmattalabi-2022} \\
 &  & \cite{Azizi-2018} \\
  \midrule
Most Vulnerable First & MVF & \cite{Kaya-2023} \\
 \midrule
 \multirow{2}{*}{Randomized Most Idle} & \multirow{2}{*}{RMI} & \cite{Mandelbaum-2012} \\
 & & \cite{Tseytlin-2009}\\
 \midrule
\multirow{5}{*}{Shortest-Queue First} & \multirow{5}{*}{SQF}& \cite{Ephremides-1980} \\
& & \cite{Johri-1989} \\
& & \cite{Hordijk-1990} \\
& & \cite{Koole-1999} \\
& & \cite{Akgun-2011} \\
\bottomrule
\label{t:literature}
\end{tabular}
\end{table}

\section{Methods}
Upon a thorough review of the existing literature, we have identified the most appropriate queue routing strategies for the housing and support service system for homeless youth in New York City (NYC). Taking into account insights obtained from our discussions with service providers, we introduce three additional novel strategies that could potentially enhance the equitable coordination of services: (i) \textit{``Greatest Number of Needs Served First"} (GNNSF) strategy, which prioritizes assigning youth to the server pool that is capable of catering to the highest number of services requested by them, (ii) \textit{``Greatest Number of Needs Served First-Idle"} (GNNSF-ID) strategy that routes youth to the server pool with an idle server that provides the greatest number of RHY's service requests, and (iii) \textit{``Largest Number of Idle Servers First"} (LNISF) strategy, where a youth is assigned to the server pool with the largest number of idle (unoccupied) servers. We summarize the abbreviated forms of these routing strategies in Table \ref{t:methods} and outline the information required to assign homeless youth to a suitable shelter.

\begin{table}[h]
\centering
\caption{The queue routing strategies we consider in this study due to their relevance. }
\begin{tabular}{lcl}
\toprule
\textbf{Routing Strategy} & \multicolumn{1}{l}{\textbf{Abbreviation}} & \textbf{Required Information Upon Arrival} \\ 
\midrule
\multirow{3}{*}{\begin{tabular}[c]{@{}l@{}}Greatest Number of Needs \\ Served First\end{tabular}} &	\multirow{3}{*}{GNNSF} & the list of services requested by each\\
& & youth and the list of services provided \\ 
& & by each shelter\\
\midrule
\multirow{3}{*}{\begin{tabular}[c]{@{}l@{}}Greatest Number of Needs \\ Served First-Idle \end{tabular}} &	\multirow{3}{*}{GNNSF-ID} & list of services requested by each youth,\\
& & list of services provided by each shelter,\\ 
& & and the idleness at shelters (idle or not)\\
\midrule
Largest Number of Idle Servers   &	LNISF & the number of idle beds at each shelter\\
\midrule
{Longest-Idle Server Pool First} & {LISF} & arrival and departure times of each RHY \\
 \midrule
%Most Vulnerable First & MVF & demographic information of each RHY \\
%\midrule
Randomized Most-Idle    &	RMI & the number of idle beds at each shelter\\
 \midrule
Shortest-Queue First & SQF & queue length at each shelter\\
\bottomrule
\label{t:methods}
\end{tabular}
\end{table}

\subsection{Routing Strategies}
% We simulate the service process of a set of RHY crisis-emergency shelters in NYC. These shelters are funded by a single agency and use coordinated entry system (CES) that enables service providers to access information regarding the available services and resources in different shelters within the system. CES also allows service providers to evaluate, and assign RHY to these available housing and support services such as mental health support, financial assistance, and substance abuse support. 

Usually when a RHY $y\in Y$ arrives to a shelter $s\in S$, they seek urgent housing and a set of support services ($R_y$) that will eventually aid in rehabilitation. On the other hand, due to limited resources, RHY often do not get to access all the services they request from the shelters. The list of services provided in shelters ($I_s$) greatly vary from shelter to shelter. Moreover, shelters often have inclusion and exclusion criteria that restrict who can or cannot access their services. Let $D$ represent the set of demographic characteristics. In this study we define the list of accepted demographic profiles for shelter $s$ as $a_s \subseteq D$ and the demographic profiles of youth as $d_y \subseteq D$. If a youth $y$ cannot receive services from any of the shelters due to demographic mismatch, they cannot be assigned to a bed. While modeling, we add them to the $L^{mismatch}$ set to keep track of the demographics that have reduced access to public housing and support services. If there are shelters in the system that accept RHY's demographic, RHY are assigned to a shelter $s \in S$, and we define this shelter as $z_y$. 

\begin{figure}[h] 
\begin{centering}
\includegraphics[height=2.8in]{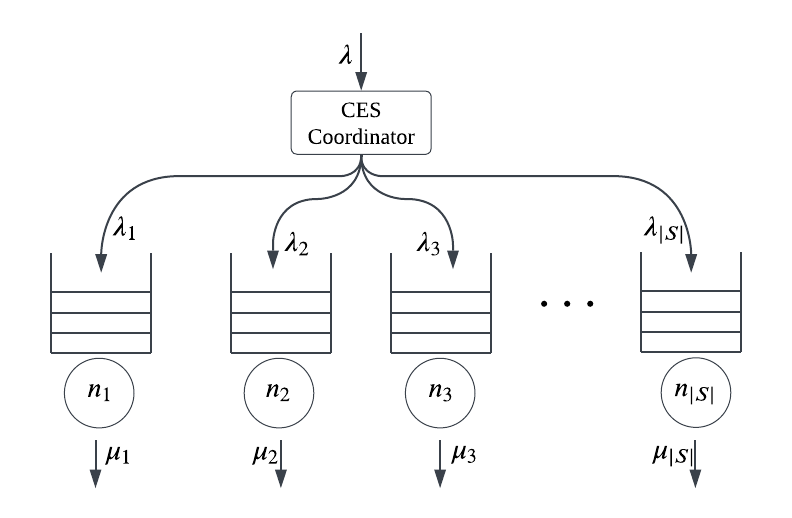}
\caption{The NYC crisis-emergency system illustrated as the inverted-V model.} 
\label{f:inverted-v-model}
\end{centering}
\end{figure}

Considering both the demographic restrictions and service availability at shelters, we simulate this system as an inverted-V queuing network with $|S|$ queues, as depicted in Figure \ref{f:inverted-v-model}. We model each of these individual queues as an $M/G/n_s+G$ queue, where $M$ represents the exponentially distributed inter-arrival times with a mean of $1/\lambda_s$, $G$ represents the generally distributed service times with mean $1/\mu_s$, $n_s$ is the number of beds in shelter $s$, and $+G$ is the generally distributed patience times of RHY with the mean $\theta$. Thus, our model considers $n_s$ homogeneous servers within $|S|$ heterogeneous server pools serving $|Y|$heterogeneous customers who have limited patience and can abandon the system before service begins.

%\vspace{-6mm}

The sets and parameters we use in this study are summarized in Tables \ref{t:sets} and \ref{t:parameters}, respectively. The queuing performance metrics we use are presented in Table \ref{t:qpm}. After modeling this housing allocation problem as a queuing network, we compare the routing strategies given in Table \ref{t:methods} to assign RHY $y$ to shelter $s$. The rest of this section discusses the routing strategies and provides pseudo-codes to illustrate our modeling process. For each routing strategy, if multiple shelters are tied (e.g., if multiple have the shortest queue length in SQF or if multiple shelters serve the maximum number of needs of a RHY), we choose the shelter to assign RHY randomly among the tied shelters. 

\begin{table}[h]
\centering
\caption{List of sets and corresponding definitions used in our study.}
\begin{tabular}{cl}
\toprule
\textbf{Symbol}  & \textbf{Definition} \\ 
\midrule
$D$ & Set of demographics\\
$Y$ & Set of RHY arriving independently to system, where $y \in Y$  \\
\multirow{2}{*}{$L^{mismatch}$} & Set of RHY that cannot be served at any shelter due to demographic mismatch, \\
& where $L^{mismatch} \subseteq Y$ \\
$L^{serviced}_s$ & Set of RHY that are served at shelter $s$, where $L^{serviced}_s \subseteq Y$ \\
\multirow{2}{*}{$L^{abandoned}_s$} & Set of RHY that abandoned shelter $s$ before service due to long wait times, \\
& where $L^{abandoned}_s \subseteq Y$ \\
$S$ & Set of shelters in the system, where $s \in S$\\
$S'_y$ & Set of shelters that can serve RHY $y$, where $S'_y  \subseteq S$\\
$I_s$ & Vector of services provided in shelter $s$  \\
$R_y$ & Vector of services requested by RHY $y$  \\
$T_{y,s}^{a} $ & List of ``$a$" times, where ``$a$" is arrival or departure before RHY $y$ arrives to shelter $s$\\
\bottomrule
\label{t:sets}
\end{tabular}
\end{table}

\begin{table}[H]
\centering
\caption{List of parameters and their definitions used in our study.}
\begin{tabular}{cl}
\toprule
\textbf{Symbol}  & \textbf{Definition} \\ 
\midrule
$n_s$ & Number of beds in shelter $s$ \\
$d_y$ & Demographic profile of RHY $y$ \\
$a_s$ & Accepted demographic profile of shelter $s$ \\
$z_y$ & The shelter RHY $y$ is assigned to\\
$t^{a}_y$ & The time of ``$a$" for RHY $y$ where ``$a$" is arrival, departure, wait, etc. time of RHY $y$\\
$t^{idle}_s$ & Idle time at shelter $s$\\
\midrule
$\mu_s$ & Mean service rate of shelter $s$ (1/length of stay)\\
$\lambda_s$ & Mean arrival rate to shelter $s$ \\
$\theta$ & Mean patience time of youth\\
$\sigma_s$ & Standard deviation of service time at shelter $s$\\
$\sigma_\theta$ & Standard deviation of patience time\\
\bottomrule
\label{t:parameters}
\end{tabular}
\end{table}

\begin{table}[H]
\centering
\caption{List of queuing performance metrics used in this study.}
\begin{tabular}{cl}
\toprule
\textbf{Symbol}  & \textbf{Definition} \\ 
\midrule
$X_{s,y}$ & Number of RHY at shelter $s$ when RHY $y$ arrives \\
$Q_{s,y}$ & Queue length at shelter $s$ when RHY $y$ arrives \\
$W_{s,y}$ & The time RHY $y$ would have to wait at shelter $s$ before receiving services \\
$\rho_{s,y}$ & Utilization of shelter $s$ when RHY $y$ arrives \\
\midrule 
$E[X_s]$ & Average number of RHY in the system at shelter $s$\\
$E[Q_s]$ & Average queue length at shelter $s$\\
$E[W_s]$ & Average wait time at shelter $s$\\
$P\{Ab_s\}$ & Proportion of RHY abandoning shelter $s$ \\
$\rho_s$ & Average long-run utilization of shelter $s$ \\
\bottomrule
\label{t:qpm}
\end{tabular}
\end{table}

\subsubsection{Greatest Number of Needs Served First:} We propose the GNNSF strategy for organizations that provide various types of services in an attempt to place youth at shelters that best match their needs. In this strategy, when a RHY $y$ arrives to the system, the coordinator first checks the set of shelters, $S'_y$, that provide housing and support services to the RHY's demographic. If there are no shelters that serve the youth's demographic, the coordinator cannot assign them to any shelters, and RHY $y$ joins the $L^{mismatch}$ set. If there is more than one shelter that accepts a RHY's demographic, the coordinator directs them to the shelter $s$ that can fulfill the most of their service needs. To do so, we compare the list of services requested by RHY $y$, $R_y$ and the list of services provided in shelter $s$, $I_s$. The pseudo-code for GNNSF is given in Algorithm \ref{alg:GNNSF}. This strategy is motivated by the appreciation that shelter and associated support services reduce RHY's vulnerability by providing security, improving their physical and mental health,  and increasing their likelihood of employment. In particular, these services disrupt trafficking activity by decreasing vulnerability for those at-risk of trafficking, including human trafficking survivors. However, this routing strategy does not account for how long a RHY may have to wait for a bed at the shelter that meets as many of its needs as possible. 

\begin{algorithm}[]
\caption{Greatest Number of Needs Served First}\label{alg:GNNSF}
\begin{algorithmic}
\State \textbf{initialize} $L^{mismatch} \leftarrow \emptyset$
\For{$ y \in Y$} 
    \If{$|S'_y | = 0$}
        \State $L^{mismatch} \leftarrow L^{mismatch} \cup y$ 
    \ElsIf{$|S'_y | > 0$}
        \State $\Delta_y \leftarrow \emptyset$ \Comment{\parbox[t]{.6\linewidth}{\textit{$\Delta_y$ stores information regarding number of needs that can be fulfilled at each shelter that accepts RHY's demographic}}}
        \For{$ s \in S'_y $}
            \State $match = I_s + R_y$ 
            \State $\Delta_y \leftarrow \Delta_y \cup |match\{i:i=2\}|$ \Comment{\parbox[t]{.5\linewidth}{\textit{Append number of requested services shelter $s$ can provide to youth $y$ }}}
        \EndFor
        \State $z_y \leftarrow argmax_s(\Delta_y)$ 
    \EndIf
\EndFor
\end{algorithmic}
\end{algorithm}

\subsubsection{Greatest Number of Needs Served First - Idle:} To address the issue of RHY potentially having to wait at a shelter that meets the greatest number of their needs while another shelter has a currently open bed, we extend the GNNSF strategy to consider idleness at shelters. When a RHY arrives to the system, the coordinator first checks which shelters accept the RHY's demographic. Of all the shelters that accept the RHY's demographic and have an idle bed available, the coordinator routes the youth to a shelter that fulfills as many of RHY's needs as possible. If there are no beds open in the system, RHY will be queued to the shelter that will eventually fulfill the greatest number of the youth's needs. With this extension, we manage to consider efficiency in addition to RHY's preferences. The pseudo-code is given in Algorithm \ref{alg:GNNSF-ID}.

\begin{algorithm}[]
\caption{Greatest Number of Needs Served First-Idle}\label{alg:GNNSF-ID}
\begin{algorithmic}
\State \textbf{initialize} $L^{mismatch} \leftarrow \emptyset$
\For{$ y \in Y$} 
    \If{$|S'_y | = 0$}
        \State $L^{mismatch} \leftarrow L^{mismatch} \cup y$ 
    \ElsIf{$|S'_y | > 0$}
        \State $\Delta_y \leftarrow \emptyset$ 
        \If{$\sum_{s \in S}(max( {n_s}-X_{s,y}, 0)) = 0$ }
            \For{$ s \in S'_y $}
                \State $match = I_s + R_y$ 
                \State $\Delta_y \leftarrow \Delta_y \cup |match\{i:i=2\}|$ 
            \EndFor
            \State $z_y \leftarrow argmax_s(\Delta_y)$
        \ElsIf{$max(\sum_{s \in S}{n_s}-\sum_{s \in S}X_{s,y}, 0) > 0$}
            \For{$ s \in \{s' \in S'_y : n_{s'} -X_{s',y} >0\}$ }
                \State $match = I_s + R_y$ 
                \State $\Delta_y \leftarrow \Delta_y \cup |match\{i:i=2\}|$ 
            \EndFor
            \State $z_y \leftarrow argmax_s(\Delta_y)$ 
        \EndIf
    \EndIf
\EndFor
\end{algorithmic}
\end{algorithm}

\subsubsection{Largest Number of Idle Servers First:}\label{ss:LNISF}  We introduce the LNISF strategy to ensure high efficiency at all shelters. This strategy increases the average number of RHY being directed to shelters with a large number of servers. We model LNISF in a way that when RHY $y$ arrives, the coordinator assigns them to the shelter $s$ that has the largest number of open beds. We present the pseudo-code for LNISF in Algorithm \ref{alg:lnisf}. 

\begin{algorithm}[]
\caption{Largest Number of Idle Servers First}\label{alg:lnisf}
\begin{algorithmic}
\State \textbf{initialize} $L^{mismatch} \leftarrow \emptyset$
\For{$ y \in Y$} 
    \If{$|S'_y | = 0$}
        \State $L^{mismatch} \leftarrow L^{mismatch} \cup y$ 
    \ElsIf{$|S'_y | > 0$}
        \State $\Delta_y \leftarrow \emptyset$ \Comment{\parbox[t]{.6\linewidth}{\textit{$\Delta_y$ stores information regarding number of idle servers at shelters that accept RHY's demographic}}}
        \For{$ s \in S'_y $}
            \State $tmp \leftarrow max(n_s-X_{s,y}, 0)$ \Comment{\parbox[t]{.55\linewidth}{\textit{$tmp$ represents the number of idle beds in shelter $s$}}}
            \State $\Delta_y \leftarrow \Delta_y \cup tmp$
        \EndFor
        \State $z_y \leftarrow argmax_s(\Delta_y)$ 
    \EndIf
\EndFor
\end{algorithmic}
\end{algorithm}

\subsubsection{Longest-Idle Server Pool First:}\label{ss:LISF}
The fourth routing strategy that we choose to assess is LISF, which is known to ensure fairness to servers \citep{Mandelbaum-2012}. In this scenario, after checking the demographic eligibility, the coordinator directs RHY to a shelter $s \in S$, where the server (bed) that has been idle the longest exists. This routing strategy, therefore, requires information regarding when each youth $y \in Y$ starts and stops receiving service at the shelter. Sometimes in systems with fairly short service times (seconds/minutes) and high customer traffic, such as emergency departments, this information might get burdensome to collect. However, in systems with longer service times (hours/days/years) such as public housing systems, information regarding each server's idle time is known to be more accessible \citep{NYC-HR-2020}. The algorithm we use to model LISF is given as Algorithm \ref{alg:lisf}.

\begin{algorithm}[]
\caption{Longest-Idle Server Pool First}\label{alg:lisf}
\begin{algorithmic}
\State \textbf{initialize} $L^{mismatch} \leftarrow \emptyset$, 
\For{$ y \in Y$} 
    \If{$|S'_y | = 0$}
        \State $L^{mismatch} \leftarrow L^{mismatch} \cup y$ 
    \ElsIf{$|S'_y | > 0$}
\State $\Delta_y \leftarrow \emptyset$ \Comment{\parbox[t]{.5\linewidth}{\textit{$\Delta_y$ stores idle time information for each shelter that accepts RHY's demographic}}}
        \For{$ s \in S'_y $}
            \State $tmp_s \leftarrow max(n_s -  X_{y,s},0)$ \Comment{\textit{$tmp_s$ represents the number of idle beds in shelter $s$}}
            \If{$tmp = 0$}
                \State $t^{idle}_s \leftarrow 0$ \Comment{\textit{If there are no open beds, idle time at shelter $s$ is $0$}}
            \ElsIf{$tmp_s > 0$}
                \State $T^{departures^*} \leftarrow \{j \in T_{y,s}^{departures}:0 < j <t^{arrival}_y  \} $
                \State $t^{idle}_s \leftarrow t^{arrival}_y -T_{y,s}^{departures \ increasing \ sorted}[|T_{y,s}^{departures}|-tmp_s]$ \Comment{\parbox[t]{.34\linewidth}{\textit{Otherwise, idle time at the shelter is the time between RHY's arrival and the longest amount of time a bed at shelter $s$ has been idle}}}
            \EndIf
                $\Delta_y \leftarrow \Delta_y \cup t^{idle}_s$
        \EndFor
        \State $z_y \leftarrow arg max_s(\Delta_y)$ 
    \EndIf
\EndFor
\end{algorithmic}
\end{algorithm}

\subsubsection{Randomized Most Idle:} \label{ss:RMI} The RMI strategy assigns RHY to an available shelter with probability equal to the fraction of idle servers in that shelter out of the total number of idle servers in the system. This strategy extends the LNISF strategy (given in Section \ref{ss:LNISF}) by incorporating randomness in assignments. In this strategy, although it is more likely, a RHY does not always end up at the shelter with the largest number of idle beds. We present the pseudo-code for this strategy in Algorithm \ref{alg:rmi}. 

% \begin{algorithm}[H]
% \caption{Randomized Most Idle}\label{alg:rmi_yaren}
% \begin{algorithmic}
% \State \textbf{initialize} $L^{mismatch} \leftarrow \emptyset$
% \For{$ y \in Y$} 
%     \If{$|S'_y | = 0$}
%         \State $L^{mismatch} \leftarrow L^{mismatch} \cup y$ 
%     \ElsIf{$|S'_y | > 0$}
%         \State $\Delta_y \leftarrow \emptyset$ \Comment{\parbox[t]{.6\linewidth}{\textit{$\Delta_y$ stores information regarding idleness ratio for shelters that accept RHY's demographic}}}
%         \State total = 0
%         \For{$ s \in S'_y $}
%             \State $tmp \leftarrow max(n_s-X_{s,y}, 0)$ \Comment{\parbox[t]{.55\linewidth}{\textit{$tmp$ represents the number of idle beds in shelter $s$}}}
%             \State $total \leftarrow total + max(n_s-X_{s,y}, 0)$
%             \State $idleness \ ratio = tmp/total$
%             \State $\Delta_y \leftarrow \Delta_y \cup idleness\ ratio$
%         \EndFor
%         \State $z_y \leftarrow Multinomial^{-1} (\Delta_y)$
%     \EndIf
% \EndFor
% \end{algorithmic}
% \end{algorithm}

% I thought it would be easier to explain my comments by writing the algorithm instead of just adding comments on the side. 
\begin{algorithm}[]
\caption{Randomized Most Idle}\label{alg:rmi}
\begin{algorithmic}
\State \textbf{initialize} $L^{mismatch} \leftarrow \emptyset$
\For{$ y \in Y$} 
    \If{$|S'_y | = 0$}
        \State $L^{mismatch} \leftarrow L^{mismatch} \cup y$ 
    \ElsIf{$|S'_y | > 0$}
        \State $\Delta_y \leftarrow \emptyset$ \Comment{\parbox[t]{.6\linewidth}{\textit{$\Delta_y$ stores information regarding idleness ratio for shelters that accept RHY's demographic}}}
        \State total = 0
        \For{$ s \in S'_y $}
            \State $tmp_s \leftarrow max(n_s-X_{s,y}, 0)$ \Comment{\parbox[t]{.55\linewidth}{\textit{$tmp$ represents the number of idle beds in shelter $s$}}}
            \State $total \leftarrow total + tmp_s$
        \EndFor
        %This allows us to get the total number of idle beds in the whole system. Then we can figure out the fraction of idle beds.
        \For{$ s \in S'_y $}
            \State $idleness \ ratio = tmp_s/total$
            \State $\Delta_y \leftarrow \Delta_y \cup idleness\ ratio$
        \EndFor
        \State $z_y \leftarrow Categorical^{-1} (\Delta_y)$
        \Comment{\parbox[t]{.55\linewidth}{\textit{Sample from the Categorical distribution}}}
    \EndIf
\EndFor
\end{algorithmic}
\end{algorithm}

\subsubsection{Shortest-Queue First:} Lastly, we include the SQF strategy. This strategy is known to improve the system efficiency by increasing the workload of the faster servers \citep{Ephremides-1980}. Therefore, in systems where servers are humans rather than objects/machines, it causes inequity among servers. We consider our servers to be the beds within the shelters, and our service time to be the length of stay of RHY. However, this is not a concern in our case because our servers are the beds within the shelters and our service time distribution is consistent among all servers (representing RHY's length of stay at the shelter). With this strategy, when youth $y$ arrives to the system, after checking the RHY's demographic eligibility, the coordinator directs RHY $y$ to the shelter $s$ with the shortest queue length $Q_{s,y}$. This strategy is presented in Algorithm \ref{alg:sqf}. 

\begin{algorithm}[]
\caption{Shortest-Queue First}\label{alg:sqf}
\begin{algorithmic}
%\Require: $\rho_s \geq 0.9, \ \forall s \in S$
\State \textbf{initialize} $L^{mismatch} \leftarrow \emptyset$
\For{$ y \in Y$} 
    \If{$|S'_y | = 0$}
        \State $L^{mismatch} \leftarrow L^{mismatch} \cup y$ 
    \ElsIf{$|S'_y | > 0$}
        \State $\Delta_y \leftarrow \emptyset$\Comment{\parbox[t]{.5\linewidth}{\textit{$\Delta_y$ stores queue length information for each shelter that accepts RHY's demographic}}}
        \For{$ s \in S'_y $}
            \State $\Delta_y \leftarrow \Delta_y \cup Q_{s,y}$ \Comment{{\textit{$Q_{s,y}$ represents the queue length at shelter $s$ when $y$ arrives}}}
        \EndFor
        \State $z_y \leftarrow arg min_s(\Delta_y)$ 
    \EndIf
\EndFor
\end{algorithmic}
\end{algorithm}

\clearpage
\subsection{Service Provision}
%Although the focus of this study is routing strategies to improve equitable access, we model the time RHY spend in the system and present our findings regarding quality-of-service (QOS) in Section \ref{s:Results}.
As part of our method of assessing a routing strategy's ability to provide equitable access, we evaluate the time RHY spend in the queue and the percentage of RHY abandoning the system before service begins. We assume that the service process of a youth after being assigned to shelter $z_y$ is the same at every shelter. After $z_y$ is determined, there are two possibilities for each youth $y$: (i) receiving service before their patience runs out during their wait to begin service, and (ii) abandoning the queue due to a long wait time and not receiving any service. If youth $y$ receives service from shelter $s$, they are added to the $L^{serviced}_s$ set, otherwise they are added to $L^{abandoned}_s$.  This process is shown in pseudo-code format in Algorithm \ref{alg:main}. 

\begin{algorithm}[]
\caption{Service process at Shelter $s \in S$ }\label{alg:main}
\begin{algorithmic}
%\For{$s \in S$}
\State \textbf{initialize} $X_{s}^{warmup} \leftarrow$ number of youth in shelter system $s$ at the end of the warmup period, 
\State $T^{y,s}{arrivals} \leftarrow \emptyset$, 
\State $T_{y,s}^{departures} \leftarrow$ departure times of youth who are in the shelter $s$ system at the end of the warmup period, 
\State $t^{arrival}_{y=0} \leftarrow 0$\\
    \For{$y \in \{y' \in Y: z_{y'}=s\}$} \Comment{{\textit{For RHY $y$ assigned to shelter $s$}}}
        \State $t^{arrival}_y  \leftarrow  t^{arrival}_{y-1}  + Exp(\lambda)$
        \State $T^{arrivals^*} \leftarrow T^{arrivals}_s$ \Comment{\parbox[t]{.4\linewidth}{\textit{Set of arrival times that occurred before $y$'s arrival}}} 
        \State $T^{y,s}{arrivals} \leftarrow T^{y,s}{arrivals} \cup t^{arrival}_y$
        %\State $T^{arrivals^*} \leftarrow \{i \in T^{arrivals}: 0 < i < t^{arrival}_y  \} $
        \State $T^{departures^*} \leftarrow \{j \in T_{y,s}^{departures}:0 < j <t^{arrival}_y  \} $   \Comment{\parbox[t]{.4\linewidth}{\textit{Set of departure times that occurred before $y$'s arrival}}}   
        \State $X_{s,y} \leftarrow X_s^{warmup}+|T^{arrivals^*}| - |T^{departures^*}|$
        \If{$X_{s,y}- n_s <0$} \Comment{\textit{If there are open beds at shelter $s$}}
            \State $t^{startservice}_y \leftarrow t^{arrival}_y$
        \Else 
            \State $\ t^{startservice}_y \leftarrow T_{y,s}^{departures \ sorted \ increasing}[|T_{y,s}^{departures}|-n_s]$
        \EndIf 
        
        \State $ t^{service}_y  \leftarrow Norm(\mu_s, \sigma_s) $
        \State $t^{patience}_y  \leftarrow Norm(\theta, \sigma_\theta)$
        \State $t^{estimatedwait}_y \leftarrow t^{startservice}_y - t^{arrival}$
        \State $\rho_{s,y} \leftarrow \frac{min(n_s, X_{s,y})}{n_s} $
        
        \If{$t^{estimatedwait}_y < t^{patience}_y $} \Comment{\textit{If youth's wait time is shorter than their patience}}
            \State $t^{actualwait}_y \leftarrow t^{estimatedwait}_y $
            \State $t^{departure}_y \leftarrow t^{startservice}_y + t^{service}_y$ \Comment{\textit{RHY leaves after receiving service}}
            \State $L^{serviced}_s \leftarrow L^{serviced}_s \cup y$
        \Else 
            \State $\ t^{actualwait}_y \leftarrow t^{patience}_y $
            \State $ t^{departure}_y \leftarrow t^{arrival}_y + t^{patience}_y$ \Comment{\textit{RHY leaves before receiving service (abandons)}}
            \State $ L^{abandoned}_s \leftarrow L^{abandoned}_s \cup y $
        \EndIf
        \State $T_{y,s}^{departures} \leftarrow T_{y,s}^{departures} \cup t^{departure}_{y-1}$
    \EndFor
%\EndFor
\end{algorithmic}
\end{algorithm}

\section{Computational Setup and Data}
We evaluate the access to four ($|S|=4$) crisis-emergency shelters that are funded by a single agency in NYC. These shelters provide a set of support services to a subset of RHY based on their demographic criteria. The demographics accepted and the set of services provided by each shelter are given in Table \ref{t:shelter_demographic} and \ref{t:shelter_services} of the Appendix, respectively. The four shelters vary in the number of beds (servers) ($n_1=53, n_2 =164, n_3 =24, n_4 =26$, where $\sum_{s\in S}n_s = 267$). 

The simulation runs for one year with an average of 2160 RHY arriving to the system each year ($|Y|=2160$) \citep{DYCD-2022}. These arriving RHY's demographic and needs profiles are created based on information gathered from various publicly available sources and meetings with stakeholders \citep{DYCD-2019,DYCD-2020,DYCD-2021}. Due to lack of consistent and comprehensive data regarding RHY's demographic profiles and needs, we could only incorporate limited information regarding their characteristics. Therefore, in this study we investigate the effect of age, gender, immigration status, and human trafficking victim status on housing access. The information we use to create the RHY demographic and needs profiles are given in Table \ref{t:youth_demographic} and \ref{t:youth_needs} of the Appendix. 

We assume that the distributions of all RHY's length of stay (i.e., service time) and the patience times are homogeneous due to lack of reliable information regarding how these parameters may differ by shelter. Following our stakeholders' recommendations, we assume the service time at shelters to be normally distributed with a mean of 60 days ($1/\mu_s =60, \forall s \in S$) and standard deviation of 5 days ($\sigma_s = 5, \forall s \in S $). The patience time of RHY follows a normal distribution with a mean of 5 days ($\theta =5$) and standard deviation of 2 days ($\sigma_{\theta} =2$). 

Meetings with stakeholders revealed that approximately 90\% of existing resources at shelters are in use by RHY on any given day. Thus, we assume that at the start of the simulation, approximately 10\% of the current resources are idle for any of the incoming youth to use. After achieving this warm-up process that fills 90\% of the capacity at all shelters, we propose a \textit{``semi-smart"} baseline strategy to represent the current CES assignment process. Using the baseline strategy, we assume that the coordinator only considers demographic eligibility and whether there is an open bed or not at shelters while placing RHY; if more than one shelter that accepts an arriving RHY's demographic has an open bed, the coordinator will assign the RHY to one of the corresponding shelters randomly. This baseline case helps us illustrate the change in queuing performance metrics compared to different routing strategies. With the baseline strategy, in total we simulate seven different routing strategies. We ran 100 instances of each of these strategies with various randomly generated youth demographic ($d_y$) and needs ($R_y$) profile combinations, considering the information provided in Table \ref{t:youth_demographic} and \ref{t:youth_needs} of the Appendix. 

\section{Results} \label{s:Results}
Here we describe the  input data used in this case study, detail the computational setup for our experiments, and present our findings regarding access to NYC RHY crisis-emergency shelters. %All experiments were conducted using Python 3.8.12 with up to 16 GB memory on an I7 processor.

\subsection{Demographic Characteristics' Effect on Access}
Many factors may reduce youth's access to public housing resources, including the systemic barriers and biases that some youth face due to their demographic characteristics \citep{Olivet-2019, Romero-2020, Fowle-2022}. \cite{HUD-2013} and \cite{EndHomelessness-2020} show that most racial and ethnic minority groups in the US experience homelessness at higher rates than white people. Another study  finds that identifying as LGBTQI+ is positively correlated with returning to homelessness within less than one year of receiving housing \citep{Petry-2021}. In addition to these disparities, our meetings with stakeholders revealed that there are very limited resources that are dedicated to RHY who are older than 21 years old and youth who have a child(ren) of their own. Motivated by these gaps in the system, to evaluate the demographic characteristics' effect on housing access we assess the overall system performance metrics for RHY from different demographic groups for each routing strategy. 

\begin{figure} []
\centering
\begin{tabular}{cc}
\includegraphics[width=0.48\textwidth]{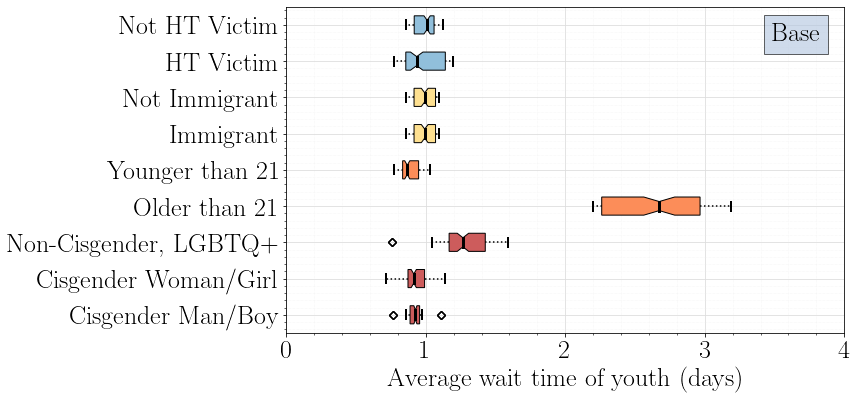} \\
\end{tabular}
\begin{tabular}{cc}
\includegraphics[width=0.48\textwidth]{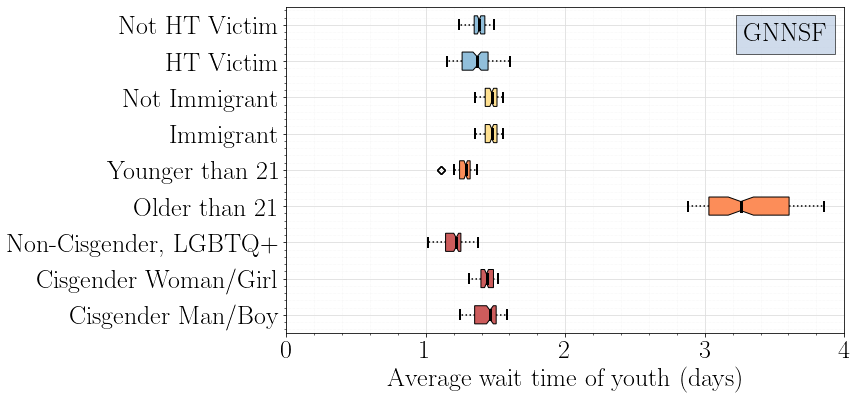}&
\includegraphics[width=0.48\textwidth]{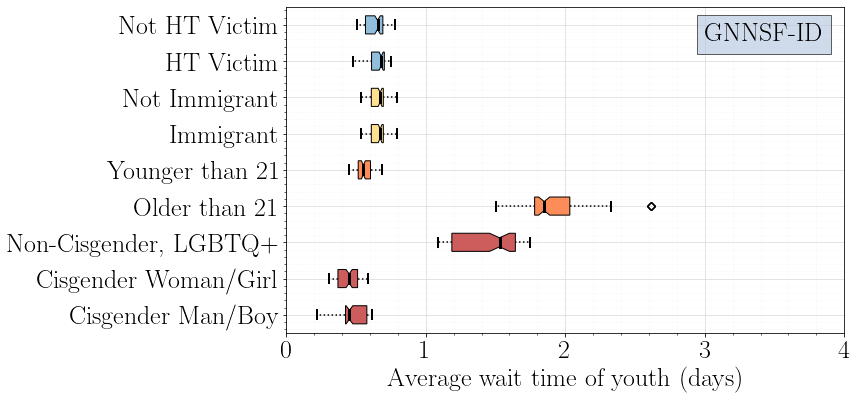}\\
\includegraphics[width=0.48\textwidth]{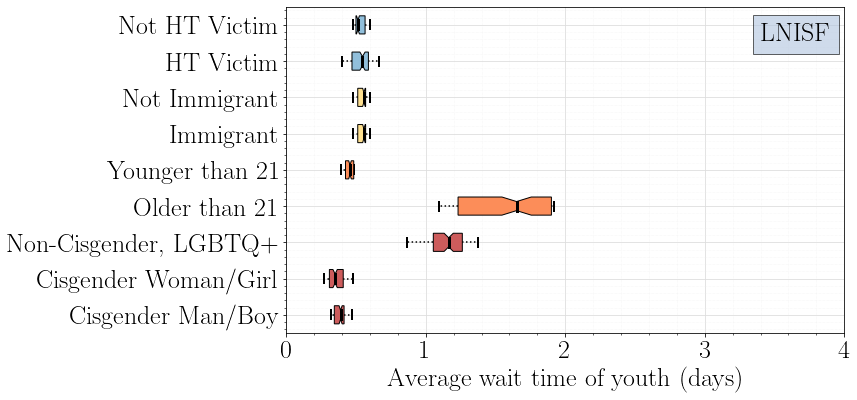} &
\includegraphics[width=0.48\textwidth]{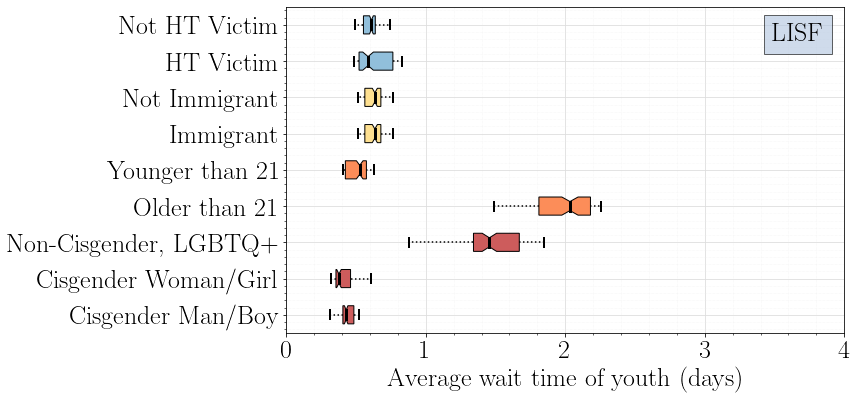} \\
\includegraphics[width=0.48\textwidth]{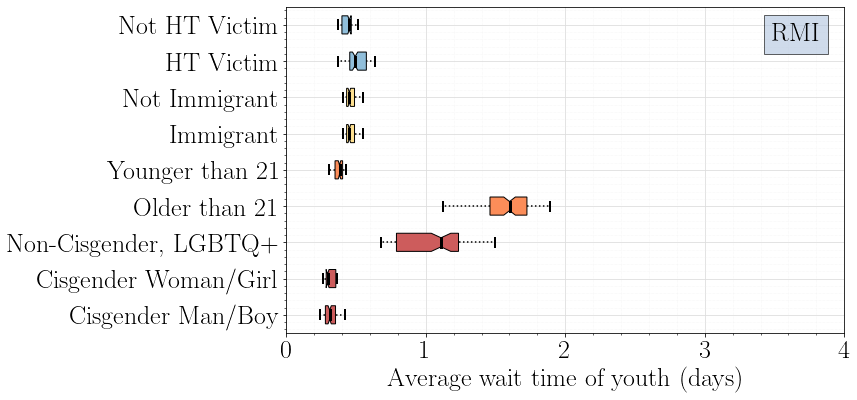} &
\includegraphics[width=0.48\textwidth]{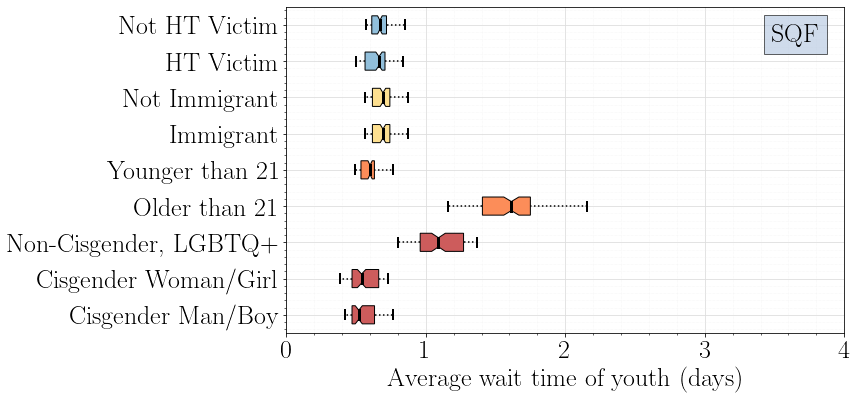} \\
\end{tabular}
\caption{The average wait time of RHY in system shown for different routing strategies and demographic characteristics considering 100 replications.}
\label{fig:demographic-wait}
\end{figure}

\begin{figure} []
\centering
\begin{tabular}{cc}
 \includegraphics[width=0.48\textwidth]{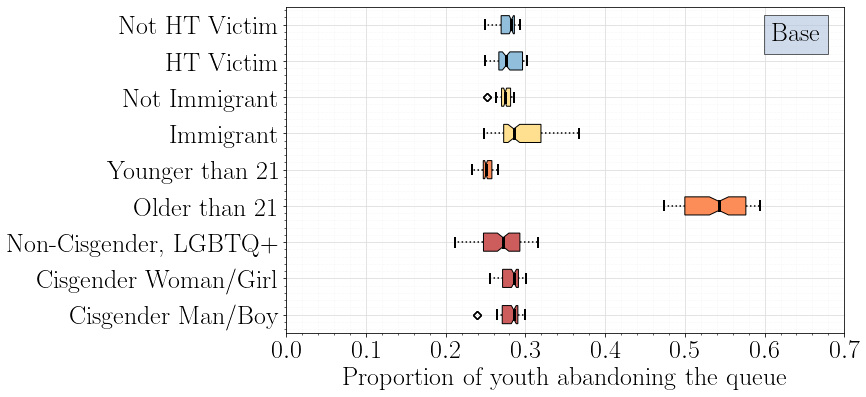}\\
\end{tabular}
\begin{tabular}{cc}
\includegraphics[width=0.48\textwidth]{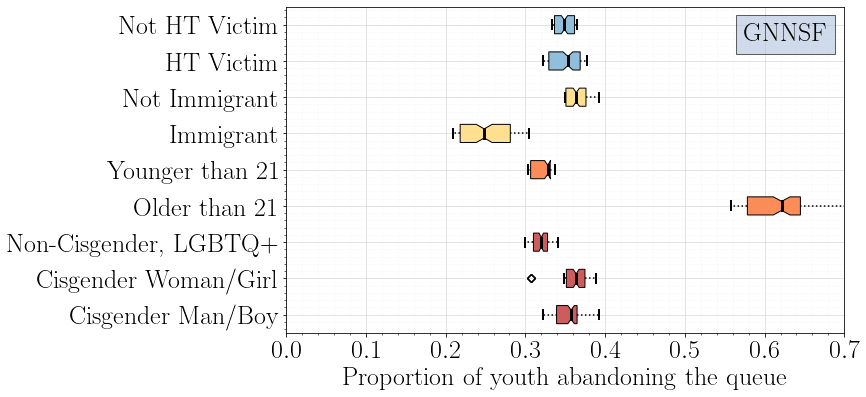}&
\includegraphics[width=0.48\textwidth]{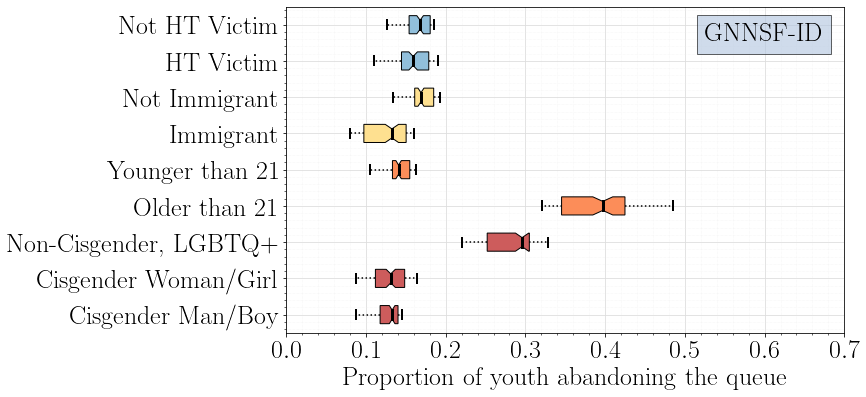}\\
\includegraphics[width=0.48\textwidth]{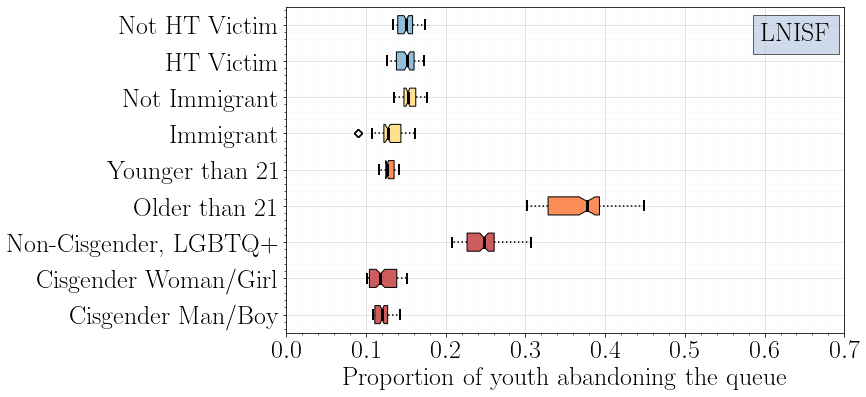} &
\includegraphics[width=0.48\textwidth]{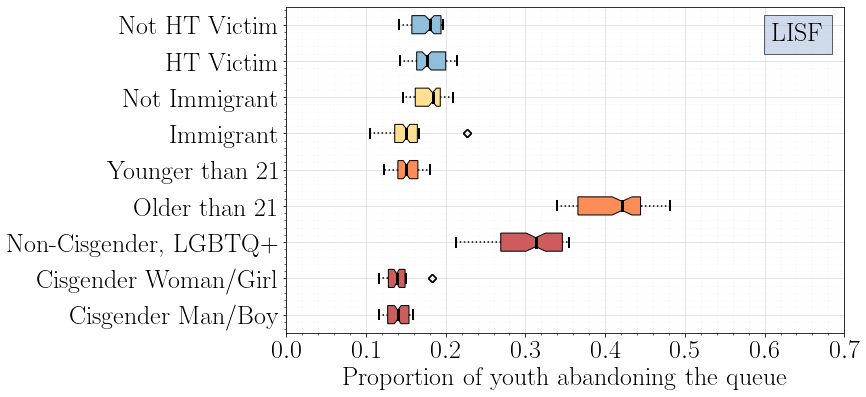} \\
\includegraphics[width=0.48\textwidth]{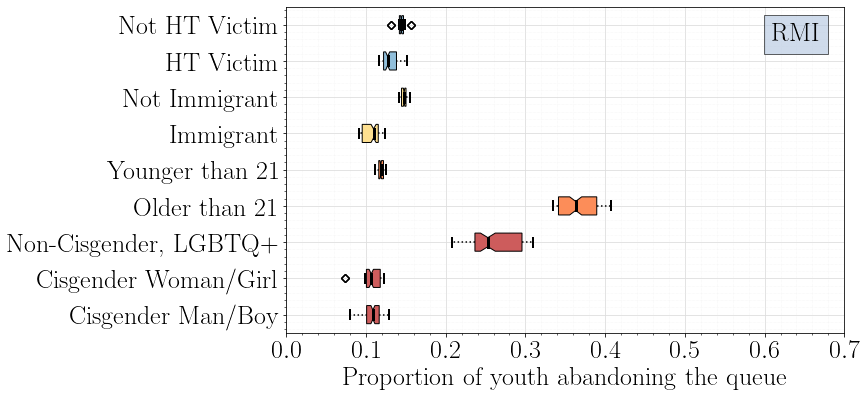} &
\includegraphics[width=0.48\textwidth]{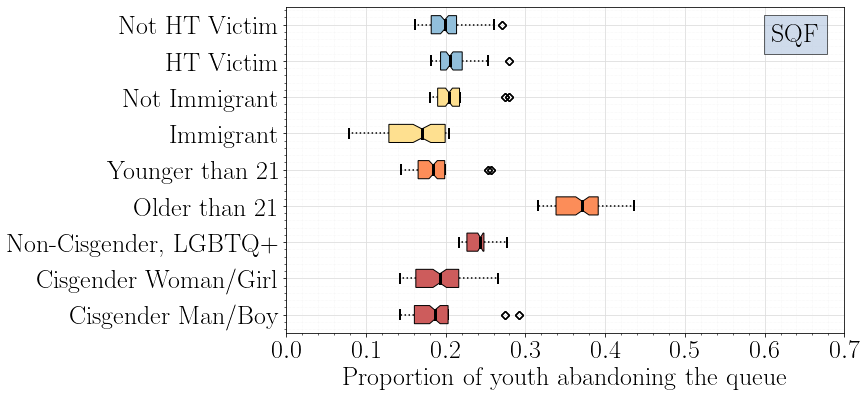} \\
\end{tabular}
\caption{The proportion of RHY abandoning the queue, shown for different routing strategies and demographic characteristics considering 100 replications.}
\label{fig:demographic-aban}
\end{figure}

Our base model estimates the average wait time of youth in the current system to be 2.07 days ($E[W] =2.07$) and the average proportion of youth abandoning as 28\% ($P\{Ab\}=28\%$), shown in Figures \ref{fig:demographic-wait} and \ref{fig:demographic-aban}. In these figures different colors of boxes represent each demographic characteristic we consider such as gender and age. %Therefore, we only compare the same colors with each other. 
The initial observation of the base scenario in each of these figures is that RHY who are older than 21 on average wait more than their younger peers and, thus, abandon the queue more. Furthermore, we see some difference in average wait times between cisgender and non-cisgender RHY. To take a closer look at this difference, we perform pairwise comparisons via t-tests and present the p-values associated with average wait time and abandonment probability comparisons in Table \ref{t:base_demographic_comparison} using a significance level of 5\% ($\alpha =0.05$). Table \ref{t:base_demographic_comparison} shows that in the current system, being older than 21 years old or identifying as non-cisgender or LGBTQ+ has a negative effect on the average wait times and abandonment proportions. 

\begin{table}[H]
\centering
\caption{Base model average wait time and abandonment proportion comparison for different demographic groups, given with their associated p-values. }
\begin{tabular}{lcc}
\toprule
Comparisons & p-value for $E[W]$  & p-value for $P\{Ab\}$\\
\midrule
Man/Boy - Woman/Girl & 0.406 & 0.970\\
Man/Boy - Non-cisgender, LGBTQ+ & 0.000 & 0.624\\
Woman/Girl - Non-cisgender, LGBTQ+ & 0.000 & 0.687\\
Older than 21 - Younger than 21 & 0.000 & 0.000\\
%With child - without child & 0.000 & 0.731\\
Immigrant - Not immigrant & 1.000 &  0.421\\
Experienced trafficking - Did not experience trafficking & 0.068 & 0.973\\
\bottomrule
\end{tabular}
\label{t:base_demographic_comparison}
\end{table}

We aim to decrease both $E[W]$ and $P\{Ab\}$ by using different routing strategies while keeping all the other parameters the same in the system. This improvement happens in different scales for the GNNSF-ID, LNISF, LISF, RMI, and SQF strategies, reducing the average wait time to 1.20, 1.15, 1.31, 0.94, and 1.54 days, respectively, while decreasing the overall abandonment proportion to 16\%, 15\%, 18\%, 14\% and 21\%, respectively.However, the GNNSF strategy fails to provide a similar improvements. In fact, it makes the metrics worse for most demographics. Although the GNNSF strategy would improve the some metrics of service quality for RHY who do receive housing services by meeting more of their stated needs, it does not perform well in terms of wait time or abandonment metrics for NYC crisis-emergency system. We hypothesize that the reason behind GNNSF's failure is the difference in number of, and types of, services provided at different shelters. 

Although these routing strategies may increase the quality of service for different demographic groups, our results indicate that the resource allocation for different demographic groups is neither balanced nor equitable to address increased vulnerabilities.  For example, our results indicate that RHY who are older than 21 have reduced access to these shelters. Moreover, if an older RHY is cisgender, there are only 26 beds that can accept their demographic (out of 267 beds total). If a cisgender, older RHY is also an immigrant, the system completely fails to provide to their demographic; thus, they wind up abandoning the system due to demographic mismatch. In each routing strategy, we observe on average 1.16\% of RHY abandoning the system simply because none of the four shelters can accommodate their demographic. The rest of this section investigates the individual effects of different demographic characteristics and provides recommendations regarding equal and equitable access. 

\subsubsection{Effect of Age}
\begin{figure} [H]
\centering
\begin{tabular}{cc}
\includegraphics[width=0.48\textwidth]{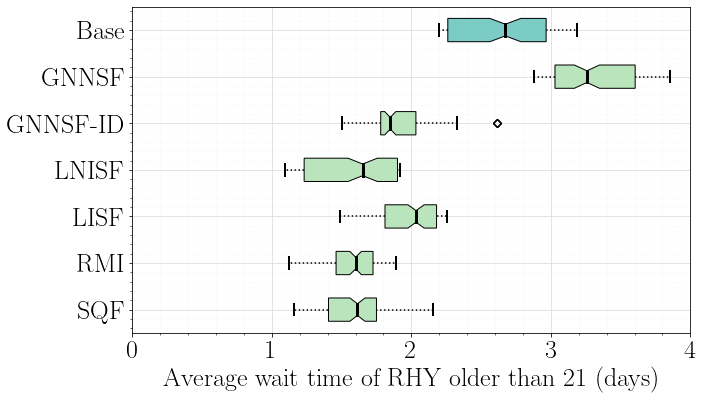}&
\includegraphics[width=0.48\textwidth]{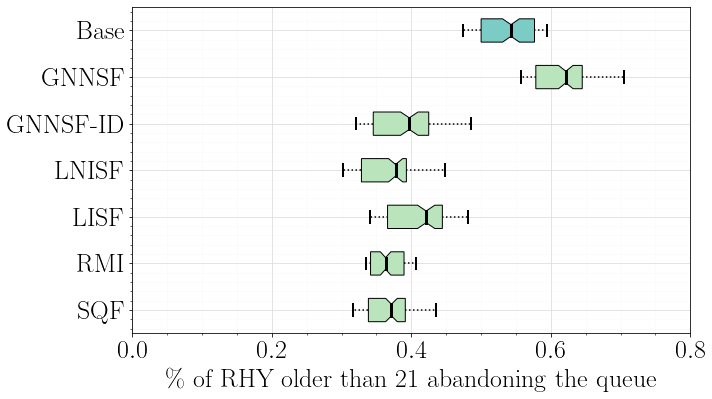} \\
\end{tabular}
\caption{The queuing performance metrics for RHY who are older than 21.}
\label{fig:age}
\end{figure}

The results of our study demonstrate that certain queue routing strategies have the potential to improve housing access for RHY who are older than 21 years old, as indicated in both Figure \ref{fig:demographic-wait} and \ref{fig:demographic-aban}. Specifically, using GNNSF-ID, LISF, LNISF, RMI, or SQF can decrease both the average wait time ($E[W]$) and the abandonment probability ($P\{Ab\}$) for this population. For example, GNNSF-ID reduces the average wait time from 2.66 days in the current system to 1.95 days, while LNISF, LISF, RMI, and SQF reduce it to 1.58, 1.96, 1.55, and 1.61 days, respectively. For ease, we display older RHY 's average wait time and abandonment proportions for different routing strategies in Figure \ref{fig:age} and display the results of our paired t-tests in Table \ref{t:age_comparison}. The p-values confirm that all the strategies except GNNSF improve the quality of service metrics for RHY who are older than 21 years old while the GNNSF strategy maintains a similar service quality to the base scenario. Our results indicate that \textbf{LNISF, RMI, and SQF} strategies improve the average wait time and the abandonment percentage the most for RHY older than 21 years old. However, even with these improvements, none of the strategies can ensure equal or equitable access for this age group due to the limited resources dedicated to them in the current system (79 beds out of 267). We argue that the resources dedicated to RHY who are older than 21 years old need to be significantly increased to achieve equality or equity for this population.

%Both in Figure \ref{fig:demographic-wait} and \ref{fig:demographic-aban} we see that using GNNSF-ID, LISF, LNISF, RMI, or SQF can improve housing access for RHY who are older than 21 by decreasing both $E[W]$ and $P\{Ab\}$. While in the current system (base scenario) the average wait time is estimated to be 2.66 days for older RHY,GNNSF-ID decreases it to 1.95 days; LNISF to 1.58; LISF to 1.96; RMI to 1.55, and SQF to 1.61 days. 
%Although they improve QOS, none of these strategies can ensure equal or equitable access for this age group due to the very limited resources dedicated for this population. In the current system there are only 79 beds that accept youth who are older than 21 years old (out of 267 beds), where 53 of them only accept LGBTQ+ RHY, and 26 of them accept any demographic group except the immigrant RHY. Therefore, to be able to talk about equality or equity for this population, we believe that the resources dedicated to them needs to be increased significantly.

\begin{table}[H]
\centering
\small
\caption{P-values for abandonment \% and average wait comparisons considering older ($>$21) RHY's access.}
\begin{tabular}{lcccccccc}
\toprule
\textbf{Comparisons} & GNNSF-Base & GNNSF-ID-Base & LNISF-Base  & LISF-Base & RMI-Base & SQF-Base \\
\midrule
\textbf{$E[W]$} & 0.000 & 0.000 & 0.000 & 0.000 & 0.000 & 0.000 \\
\textbf{$P\{Ab\}$ } & 0.098 & 0.004 & 0.001 & 0.011 & 0.001 & 0.001\\
\bottomrule
\end{tabular}
\label{t:age_comparison}
\end{table}

%  &   LISF-LNISF &LISF-SQF & LNISF-SQF \\
% \midrule
% \textbf{p-value} & 0.001 &  & 0.001 & 0.011 & 0.376 & 0.380 & 0.994 \\

\subsubsection{Effect of Gender and Sexual Orientation}
\begin{figure} [H]
\centering
\begin{tabular}{cc}
\includegraphics[width=0.48\textwidth]{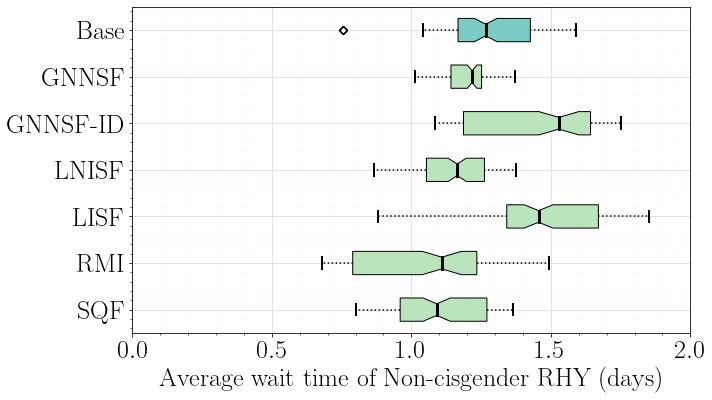}&
\includegraphics[width=0.48\textwidth]{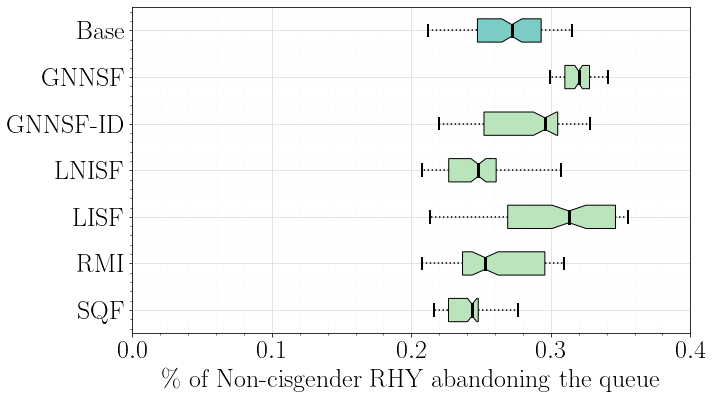} \\
\end{tabular}
\caption{The queuing performance metrics for RHY who are non-cisgender and or LGBTQ+.}
\label{fig:gender}
\end{figure}

We present our findings on the differences in housing access between cisgender and non-cisgender RHY, as well as the effectiveness of different routing algorithms in improving access. Our analysis shows that in the current system non-cisgender and LGBTQ+ RHY experience a 36\% longer average wait time compared to their straight, cisgender counterparts (1.25 days rather than 0.92), which increases their vulnerability and leads to higher abandonment rates. In Table \ref{t:gender_comparison} and Figure \ref{fig:gender} we show that the average wait time differs from the base case scenario for almost all the strategies (GNNSF's p-value is $0.056\cong0.050$). While GNNSF-ID and LISF are the only strategies that fail to  decrease the average wait time for non-cisgender youth, GNNSF, LNISF, RMI, and SQF decrease it to 1.20, 1.15, 1.05, and 1.01 days, respectively. We note that the seemingly improved access under GNNSF strategy is actually due to the high abandonment proportion of non-cisgender RHY. Therefore, we recommend employing \textbf{LNSIF, RMI, or SQF} strategies to address the access disparities and decrease the vulnerability of non-cisgender and LGBTQ+ RHY to trafficking. Our results emphasize the need to consider demographic characteristics and ensure equitable access to housing resources for all RHY.

%However, our simulations demonstrate that using LNISF, RMI, or SQF can decrease the average wait time for non-cisgender youth to 1.15, 1.05, and 1.01 days, respectively, while GNNSF-ID and LISF fail to improve access for this demographic. We note that the seemingly improved access under GNNSF is actually due to a high abandonment rate. Therefore, we recommend employing LNSIF, RMI, or SQF strategies to address the access disparities and decrease the vulnerability of non-cisgender RHY to trafficking. Our results emphasize the need to consider demographic characteristics to ensure equitable access to housing resources for all RHY.

%We display $E[W]$ and $P\{Ab\}$ for non-cisgender RHY in Figure \ref{fig:gender} to illustrate the differences between routing algorithms. We estimate that, in the current system non-cisgender RHY on average wait 36\% more than their counterparts (1.25 days rather than 0.92). Thus, more non-cisgender youth abandon their queues before receiving housing services, resulting in increased vulnerabilities. 
%When we focus only on the average wait time, the results entail that GNNSF, LNISF, RMI, and SQF can increase access to housing resources. However, shifting the focus to abandonment proportion clarifies that the wait time improvement of GNNSF is actually due to the high abandonment proportion of non-cisgender RHY. Therefore, to address these access disparities and decrease non-cisgender RHY's increased vulnerability to trafficking, we recommend routing RHY using \textbf{LNSIF, RMI, or SQF} strategies.

\begin{table}[H]
\centering
\small
\caption{P-values for abandonment \% and average wait comparisons considering non-cisgender RHY's access.}
\begin{tabular}{lcccccccc}
\toprule
\textbf{Comparisons} & GNNSF-Base & GNNSF-ID-Base & LNISF-Base  & LISF-Base & RMI-Base & SQF-Base \\
\midrule
\textbf{$E[W]$} & 0.056 & 0.000 & 0.000 & 0.000 & 0.000 & 0.000 \\
\textbf{$P\{Ab\}$} & 0.079 & 0.367 & 0.558 & 0.260 & 0.787 & 0.372 \\
\bottomrule
\end{tabular}
\label{t:gender_comparison}
\end{table}

\subsubsection{Effect of Being an Immigrant}
\begin{figure} [H]
\centering
\begin{tabular}{cc}
\includegraphics[width=0.48\textwidth]{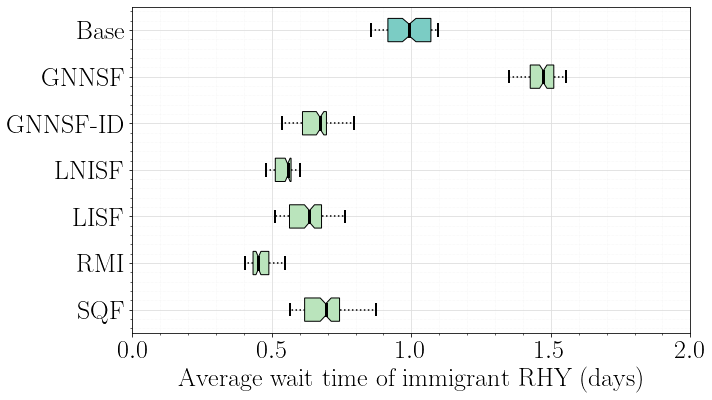}&
\includegraphics[width=0.48\textwidth]{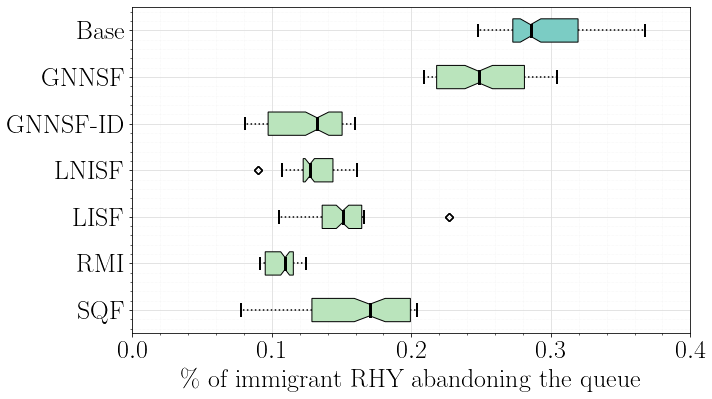} \\
\end{tabular}
\caption{The queuing performance metrics for RHY who are immigrant.}
\label{fig:immigrant}
\end{figure}

Based on the results presented in Figure 2, there is no statistically significant difference in the average wait times for immigrant RHY when compared to their non-immigrant counterparts. However, the findings presented in Figure \ref{fig:demographic-aban} raise concerns regarding the equity of abandonment proportions. Therefore, similar to our discussion on other characteristics effects, we provide an analysis of the average wait times and abandonment proportions for immigrant RHY in Figure \ref{fig:immigrant}. Our pairwise comparison of different routing strategies based on average wait time demonstrates that, unsurprisingly, all strategies offer varying levels of service quality. Specifically, GNNSF-ID and SQF reduce the average wait time for immigrant youth by an average of 30\%, LISF reduces it by 40\%, and LNISF and RMI reduce it by approximately 50\%. These findings suggest that these routing strategies can improve access to housing resources for immigrant RHY, and lead to \textit{``equitable"} access rather than equal access. Moreover, Table \ref{t:immigrant_comparison} presents the pairwise comparison of abandonment proportions and Figure \ref{fig:immigrant} confirms that \textbf{LNISF and RMI} strategies offer the most significant improvements in service quality for immigrant RHY in terms of abandonment proportions.

\begin{table}[H]
\centering
\small
\caption{P-values for abandonment \% and average wait comparisons considering immigrant RHY's access.}
\begin{tabular}{lcccccc}
\toprule
\textbf{Comparisons} & GNNSF-Base & GNNSF-ID-Base & LNISF-Base  & LISF-Base & RMI-Base & SQF-Base \\
\midrule
\textbf{$E[W]$} & 0.000 & 0.000 & 0.000 & 0.000 & 0.000 & 0.000 \\
\textbf{$P\{Ab\}$} & 0.201 & 0.000 & 0.000 & 0.000 & 0.000 & 0.000 \\
\bottomrule
\end{tabular}
\label{t:immigrant_comparison}
\end{table}

\subsubsection{Effect of Being a Human Trafficking Victim}
Table \ref{t:base_demographic_comparison} reveals that the status of being a human trafficking victim does not significantly affect youth's ability to access crisis and emergency shelters and suggests equal access for all in the base model. This equality in average wait times among human trafficking victims and their non-victim counterparts persists when we route youth using other strategies. Although these strategies are unable to provide equitable access for youth who have experienced trafficking, \textbf{all of them, except GNNSF}, improve the efficiency of the system by decreasing $E[W]$ and $P\{Ab\}$. This enhancement in the overall service quality of the system is expected to reduce the risk of vulnerable youth being trafficked or exploited. However, given the increased risk for this population to experience trafficking and exploitation, we emphasize the significance of supporting and sustaining improved, equitable access for vulnerable populations.

\begin{figure} [H]
\centering
\begin{tabular}{cc}
\includegraphics[width=0.48\textwidth]{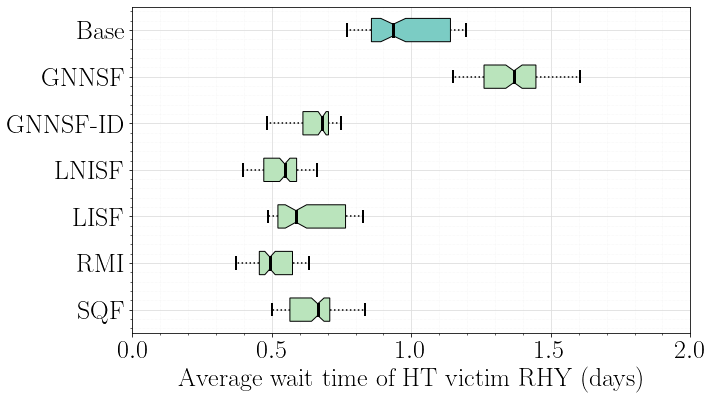}&
\includegraphics[width=0.48\textwidth]{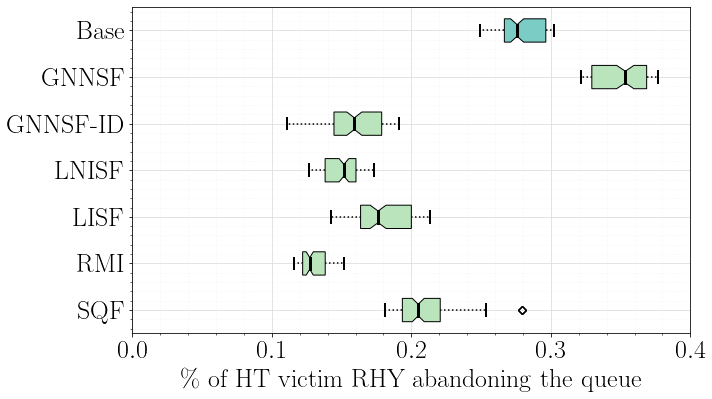} \\
\end{tabular}
\caption{The queuing performance metrics for RHY who are human trafficking victims.}
\label{fig:htvictim}
\end{figure}

\subsection{Access to Different Shelters and Fairness}
\begin{figure} []
\centering
\begin{tabular}{ccc}
\includegraphics[width=0.31\textwidth]{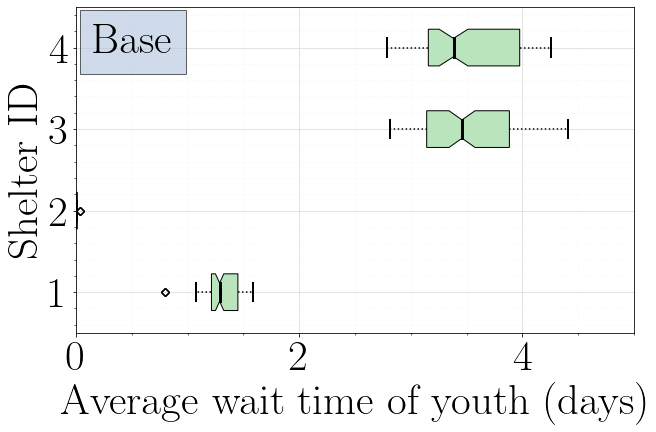}
\\
\end{tabular}
\begin{tabular}{ccc}
\includegraphics[width=0.31\textwidth]{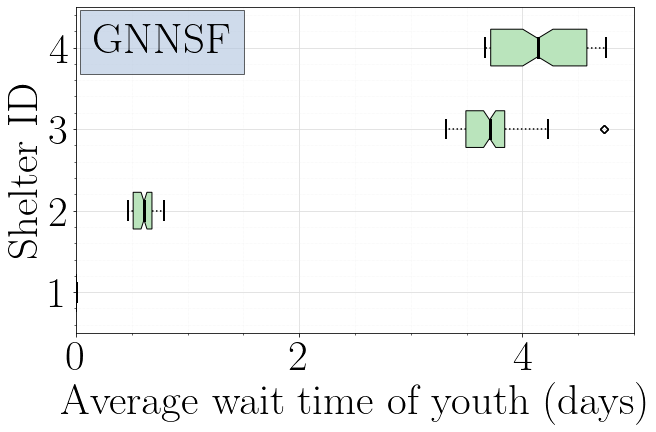} &
\includegraphics[width=0.31\textwidth]{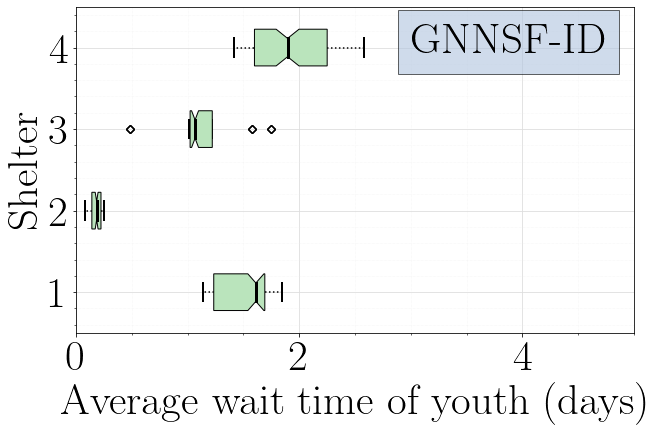} &
\includegraphics[width=0.31\textwidth]{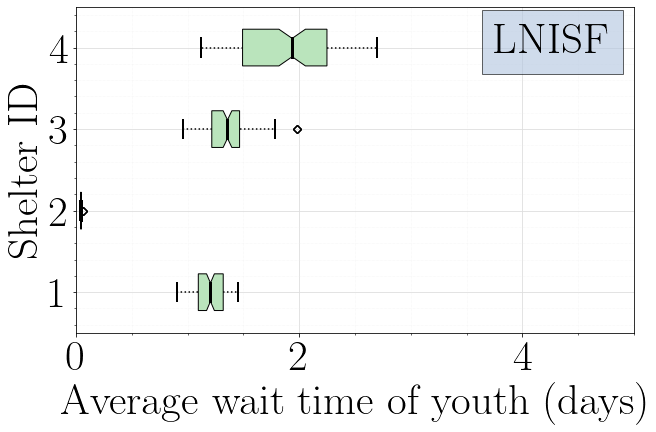} \\
\includegraphics[width=0.31\textwidth]{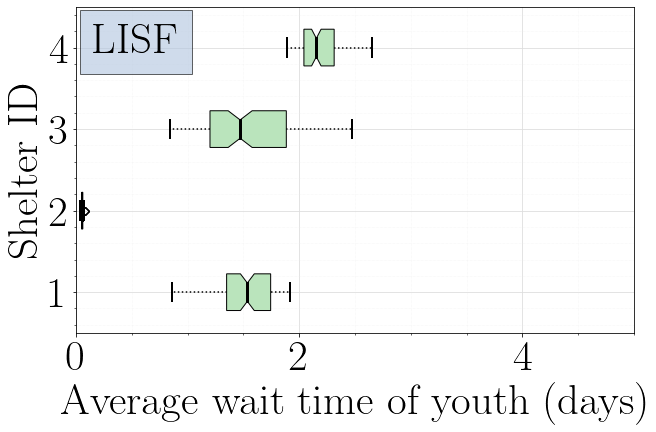} &
\includegraphics[width=0.31\textwidth]{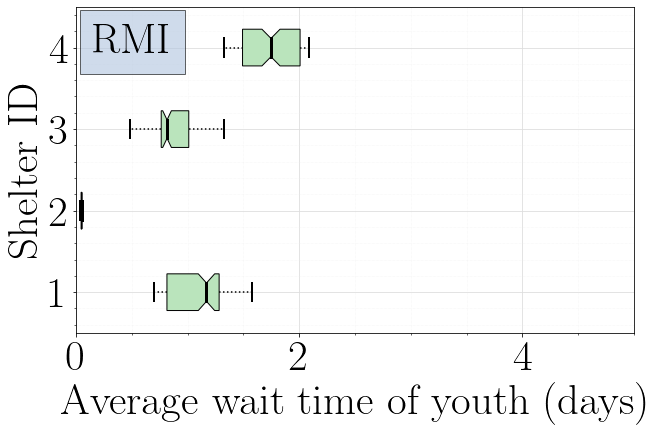}&
\includegraphics[width=0.31\textwidth]{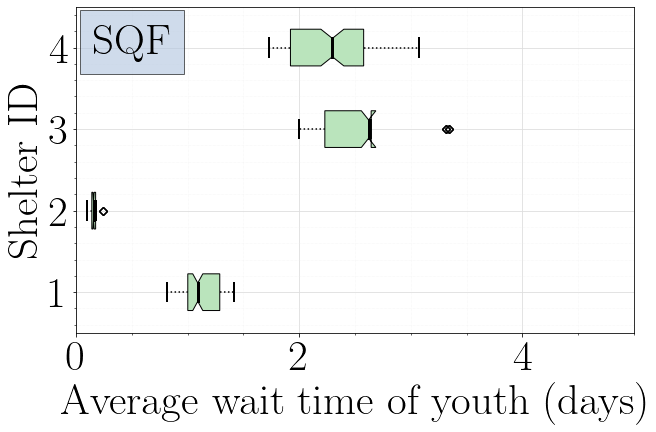} \\
\end{tabular}
\caption{The average wait time of RHY in system, shown for different routing strategies and demographic characteristics considering 100 replications.}
\label{fig:shelter-wait}
\end{figure}

\begin{figure} []
\centering
\begin{tabular}{ccc}
\includegraphics[width=0.31\textwidth]{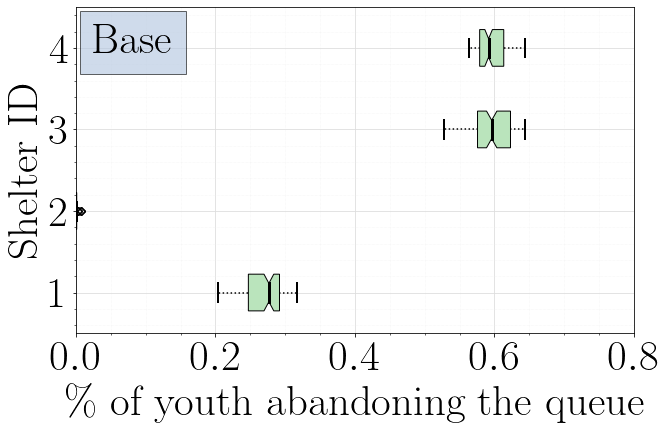}
\\
\end{tabular}
\begin{tabular}{ccc}
\includegraphics[width=0.31\textwidth]{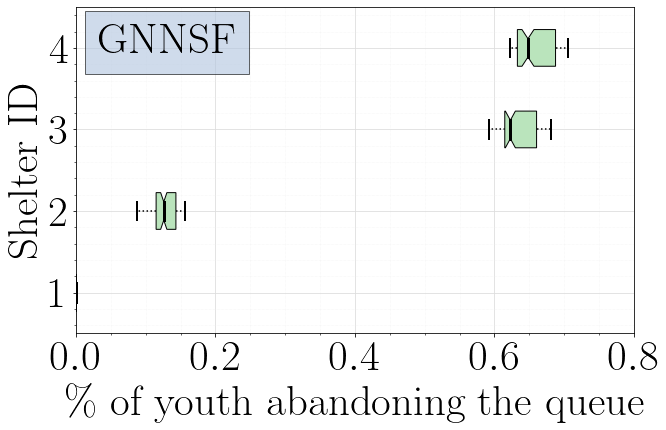} &
\includegraphics[width=0.31\textwidth]{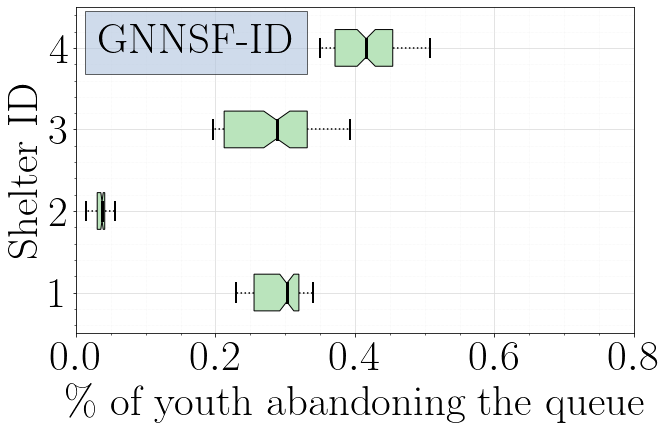} &
\includegraphics[width=0.31\textwidth]{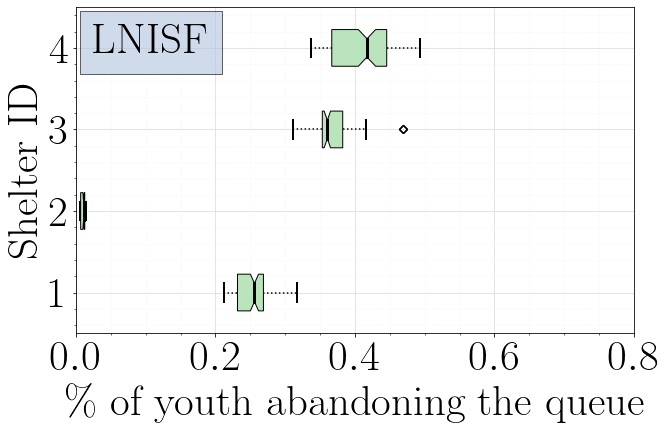} \\
\includegraphics[width=0.31\textwidth]{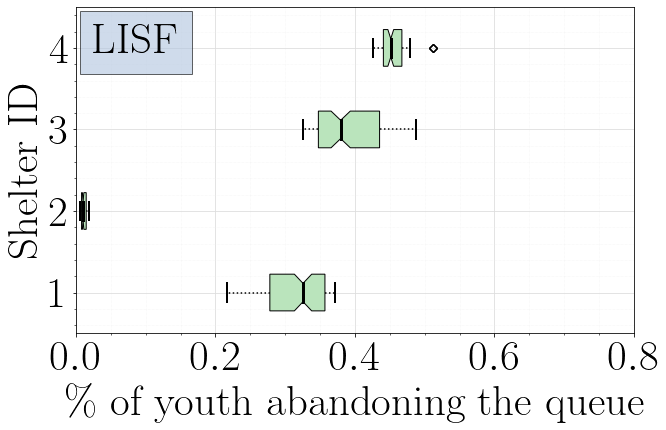} &
\includegraphics[width=0.31\textwidth]{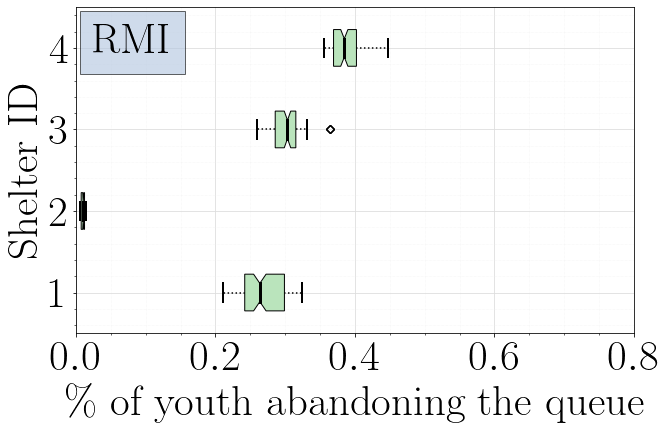}&
\includegraphics[width=0.31\textwidth]{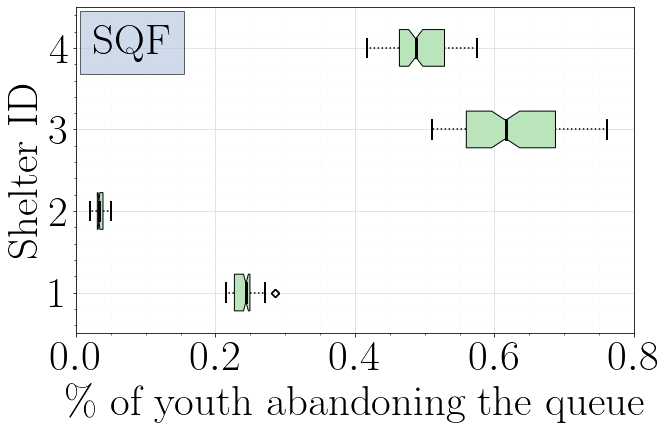} \\
\end{tabular}
\caption{The average proportion of RHY abandoning the system, shown for different routing strategies and shelters, considering 100 replications.}
\label{fig:shelter-aban}
\end{figure}

\begin{figure} []
\centering
\begin{tabular}{ccc}
\includegraphics[width=0.31\textwidth]{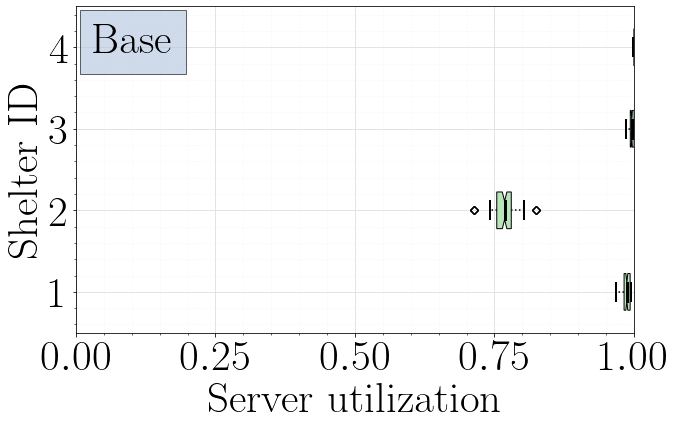}
\\
\end{tabular}
\begin{tabular}{ccc}
\includegraphics[width=0.31\textwidth]{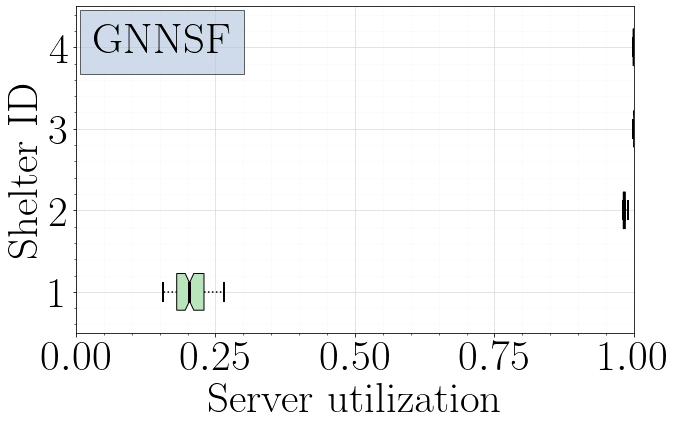} &
\includegraphics[width=0.31\textwidth]{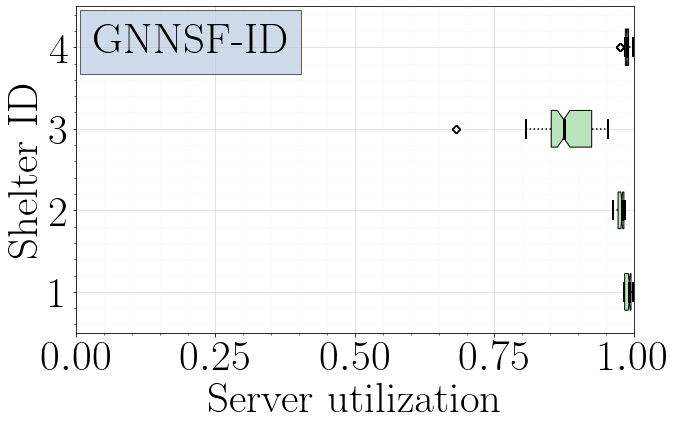} &
\includegraphics[width=0.31\textwidth]{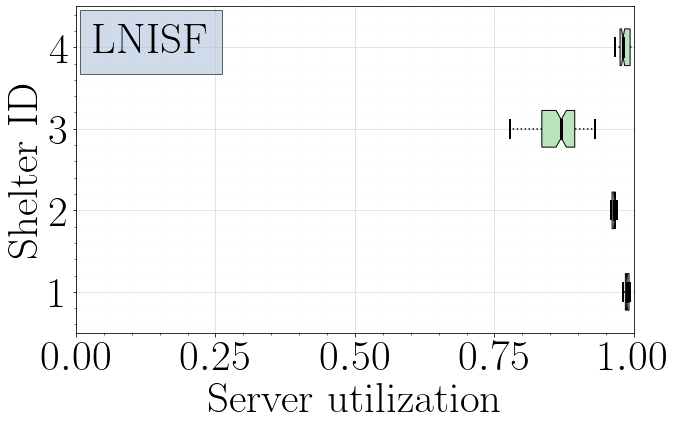} \\
\includegraphics[width=0.31\textwidth]{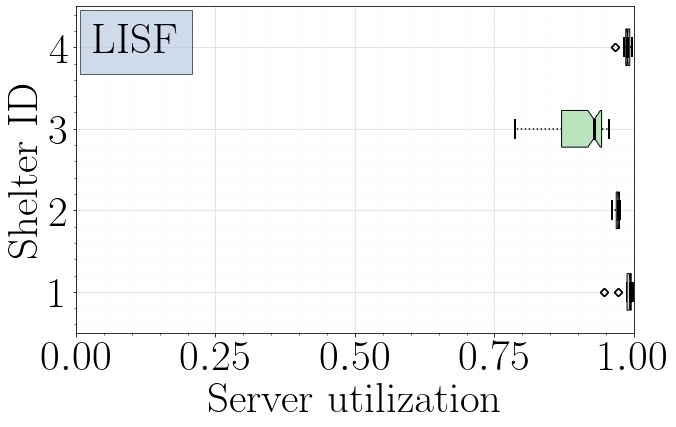} &
\includegraphics[width=0.31\textwidth]{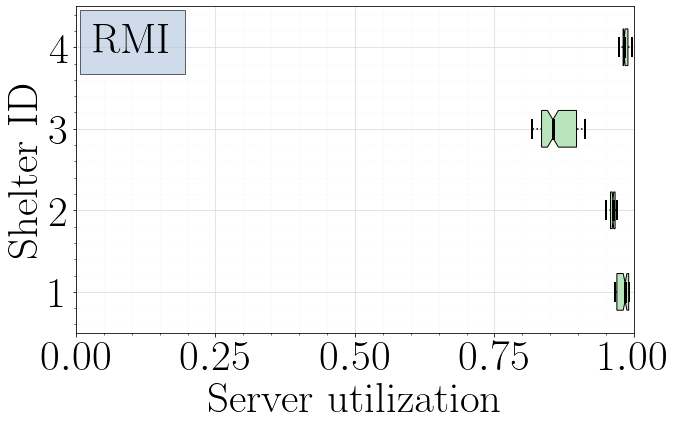}&
\includegraphics[width=0.31\textwidth]{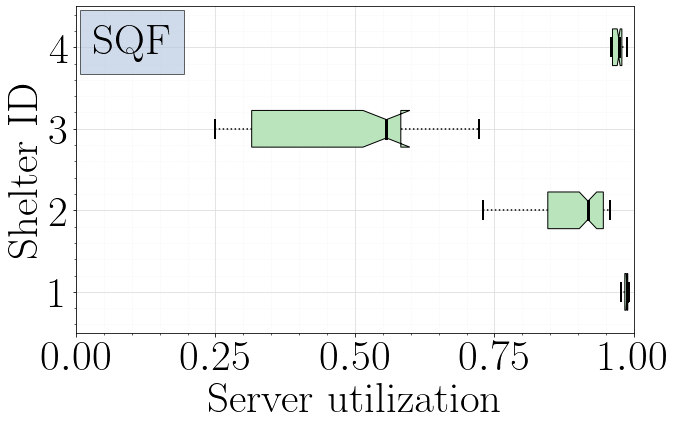} \\
\end{tabular}
\caption{The average utilization at shelters, shown for different routing strategies, considering 100 replications.}
\label{fig:shelter-utilization}
\end{figure}

In addition to examining equity and equality from RHY's perspective, we also examine fairness from the perspective of service providers, with the goal of achieving standardized service quality at various shelters. By doing so, we aim to ensure staff responsible for the same tasks at different shelters to have a similar workload. To assess this, we compare the queuing performance metrics of different shelters under various routing strategies, as depicted in Figures \ref{fig:shelter-wait}, \ref{fig:shelter-aban}, and \ref{fig:shelter-utilization}. These figures respectively present the average wait time, average proportion of abandonment, and average utilization at different shelters. Through this analysis, we can identify disparities in service quality across different shelters and identify strategies to promote fairness and equity in the delivery of services to RHY.

In our base scenario we estimate average wait time at different shelters to be 1.3, 0.0, 3.5, and 3.5 days at shelters 1 -4, respectively. This trend is also observed in the proportion of RHY abandoning and the server utilization metrics, highlighting the service quality differences among the shelters. Our findings suggest that this difference in service quality in our base case scenario is mainly driven by the number of beds available at each shelter and the demographic eligibility criteria of shelters. These results are particularly relevant for the staff at smaller shelters, who may experience a disproportionate burden due to unbalanced workloads. Therefore, we believe that our base scenario effectively captures the disparities in the system and can inform efforts to improve service quality across different shelters.

%Therefore, we believe that our base case scenario achieves to reflect the inefficiencies of the system and the burden experienced by staff at smaller shelters due to unbalanced workloads.
In our analysis of demographic characteristics’ impact on access, we observe that the LNISF and RMI routing strategies offer the highest service quality across all demographics. Figures \ref{fig:shelter-wait} and \ref{fig:shelter-aban} illustrate that these strategies consistently deliver average wait times of less than 2 days and abandonment rates of less than 50\% at each shelter. Moreover, while comparing service quality levels across different shelters, we note that although not as high as the LNISF and RMI strategies, the LISF strategy maintains a high level of service quality at each shelter. Furthermore, Figure \ref{fig:shelter-utilization} examines system efficiency and indicates that the average utilization at different shelters while using the LNISF, LISF or RMI routing strategies is also approximately the same (around 99\% at shelters 1, 2, and 4; and 90\% at shelter 2).This means balanced work load for staff working at different shelters and ensures fairness for servers.

% \subsection{Survey?}
% \begin{itemize}
% \item Do you think coordinated entry system ensures equitable and fair access to all youth?
% \item Do you think coordinated entry system is efficient?
% \item Do you have recommendations to improve coordinated entry system?
% \item What do you think about fairness and efficiency metrics of these different routing strategies that we present?
% \end{itemize}

\section{Conclusions}
%summary
%We model the operational challenge of assigning runaway and homeless youth (RHY) to scarce public housing and support resources. We use a simulation-based technique to evaluate the system efficiency, equitable access, and fairness towards servers while using different queue routing strategies. To the best of our knowledge, our study is the first attempt in the operations research and analytics literature to compare different routing strategies for a better (more efficient and equitable) housing allocation process. 

Our study aims to address the challenge of allocating scarce public housing and support resources to runaway and homeless youth (RHY) in an efficient and equitable manner. To achieve this, we employ a simulation-based approach to evaluate the effectiveness and fairness of various interpretible queue routing strategies. We consider queue routing strategies that have already been introduced to literature, while also introducing three new strategies.  To the best of our knowledge, our research is the first attempt in the operations research and analytics literature to compare different interpretible routing strategies for improving the housing allocation process in terms of efficiency and equity.

In this study, we present our simulation models that involve a queuing network with pools of multiple parallel servers. Our models account for servers with demographic eligibility criteria; stochastic youth arrival; impatient youth behaviour (possibility of abandonment); and a decision-maker (coordinator) that determines which server pool youth is routed to. In each simulation model the coordinator routes RHY differently, and according to: \textit{``Greatest Number of Needs Served First"} (GNNSF), \textit{``Greatest Number of Needs Served First - Idle"} (GNNSF-ID), \textit{``Largest Number of Idle Servers First"} (LNISF), \textit{``Longest-Idle Server Pool First"} (LISF), \textit{``Randomized Most Idle"} (RMI), and \textit{``Shortest-Queue First"} (SQF). While we illustrate our approach with a case study aimed at improving equitable access to crisis-emergency housing resources in New York City (NYC), the same (or a similar) model could be readily used for non-profit shelter organizations in other locations. Furthermore, similar simulation models could be employed by governmental and non-profit organizations that provide public services such as mental health, and substance abuse support to determine more efficient and equitable ways to assign individuals to scarce resources.

%Findings
We discovered new insights that have the potential to impact the efforts to achieve equitable housing access in NYC. Our findings indicate that there is a need for improvement in the NYC crisis-emergency system's capacity allocation and assignment strategies. The inclusion and exclusion criteria (demographic eligibility criteria) used at shelters often make it challenging for RHY to access housing and certain support services such as childcare and human trafficking (HT) related rehabilitative services. In the current system, there are very few resources that are dedicated to RHY who are LGBTQI+ and older than 21 years old. This situation increases the access disparities among different demographic groups; thus, should be addressed by increasing the housing capacity dedicated to these populations and by re-evaluating eligibility criteria at shelters that do not serve these demographics. 

In this study, we aim to increase the accessibility to scarce crisis-emergency resources and improve system efficiency only by changing the way RHY are routed to the shelters. Compared to the base scenario that represents the current NYC system, each of the routing strategies we consider provide different levels of service quality to RHY. The \textit{``Largest Number of Idle Servers First"} and \textit{``Randomized Most Idle"} strategies consistently perform the best by decreasing the average wait time of RHY by approximately one day and the overall abandonment rate of the system by 14\%. Furthermore, although the shelters in our system are very different from each other (demographics they serve, the trainings their staff go through, size of the shelters, services they provide and more), these strategies ensure balanced work load among different shelters, while increasing the system efficiency (the average server utilization at every shelter is more than 90\%). Thus, it is most likely that LNISF and RMI strategies will increase system efficiency while maintaining fairness for servers and high quality of service for RHY when implemented. Although to a lesser extent, the \textit{``Shortest-Queue First"} and the \textit{``Longest-Idle Server Pool First"} strategies also provide a higher service quality compared to the base case. However, SQF fails to maintain fairness among server pools (shelters) and the LISF strategy increases the average wait time and abandonment proportion of non-cisgender and LGBTQI+ RHY. 

Out of the remaining two routing strategies, only the \textit{``Greatest Number of Needs Served First"} strategy does not produce an increase in the service quality. Nonetheless, we believe that the incapability of GNNSF strategy can be explained by the major difference in number and type of services provided at different shelters. Thus, when every shelter in the system is capable of providing an adequate level of support services to all RHY, the GNNSF strategy may be able to improve access to RHY while also increasing their likelihood of exiting homelessness. 

%Future directions
While this study specifically focuses on routing strategies to improve equitable access to RHY in NYC, there remains great potential for organizations in other locations if the sufficient information regarding the system system parameters exist. Further extensions to this study include a queue routing strategy that considers equitable assignment based on the \textit{``housing outcomes"} rather than the assignment itself. Moreover, there remains great potential to embed our simulation models into a decision-support tool to further facilitate the decision-making process in public service settings.

%Concluding remarks
Our approach benefits government and nonprofit decision-makers by offering a means to effectively evaluate the allocation of scarce resources. Overall, our study represents an innovative use of simulation modeling to address a resource assignment problem with a broader societal impact.

\ACKNOWLEDGMENT{The authors would like to thank the New York City Mayor's Office and New York Coalition for Homeless Youth for their insights. Special thanks to Dr. Andrew Trapp, Dr. Renata Konrad, and Geri Dimas for their time and efforts. This work is supported by the National Science Foundation under Grant No. CMMI-1935602.}

\bibliographystyle{informs2014}

\bibliography{references}
\clearpage

\begin{APPENDICES}
\huge \textbf{Appendix} \\ \normalsize

\begin{table}[H]
\centering
\caption{The demographic groups accepted at each shelter.}
\begin{tabular}{lcccc}
\toprule
\textbf{Demographic Characteristics} & \textbf{Shelter 1} & \textbf{Shelter 2} & \textbf{Shelter 3} & \textbf{Shelter 4} \\
\midrule
Age Limit & 24 & 21 & 21 & 24 \\
Cisgender Man/Boy & 0 & 1 & 1 & 1 \\
Cisgender Woman/Girl & 0 & 1 & 1 & 1 \\
Transgender Man/Boy & 1 & 1 & 1 & 1 \\
Transgender Woman/Girl & 1 & 1 & 1 & 1 \\
Non-binary & 1 & 1 & 1 & 1 \\
Genderqueer & 1 & 1 & 1 & 1 \\
Immigrant Status & 1 & 1 & 1 & 0 \\
Human Trafficking Victim Status & 1 & 1 & 1 & 1\\
\bottomrule
\end{tabular}
\label{t:shelter_demographic}
\end{table}

\begin{table}[H]
\centering
\caption{The list of services provided at each shelter.}
\begin{tabular}{lcccc}
\toprule
\textbf{Support Services} & \textbf{Shelter 1} & \textbf{Shelter 2} & \textbf{Shelter 3} & \textbf{Shelter 4} \\
\midrule
Mental Health Support & 0 & 1 & 1 & 1 \\
Medical Support & 0 & 0 & 0 & 0 \\
Substance Abuse Support & 0 & 1 & 1 & 1 \\
Crisis 24-Hour Services & 1 & 1 & 1 & 1 \\
Long-Term Housing Support & 1 & 1 & 1 & 1 \\
Legal Assistance & 0 & 1 & 1 & 0 \\
Service Coordination & 1 & 1 & 1 & 1 \\
Practical Assistance & 1 & 1 & 1 & 1 \\
Financial Assistance & 1 & 1 & 0 & 1 \\
Life Skills Support & 1 & 1 & 1 & 1 \\
Employment Assistance & 0 & 1 & 1 & 1 \\
Education Assistance & 0 & 1 & 1 & 1 \\
Childcare Support & 0 & 1 & 1 & 0\\
\bottomrule
\end{tabular}
\label{t:shelter_services}
\end{table}

\begin{table}[]
\centering
\caption{The proportion of RHY who identify themselves as a certain demographic group. }
\begin{tabular}{lcc}
\toprule
Demographic Attributes & Options & \multicolumn{1}{l}{Proportion of RHY} \\
\midrule
\multirow{9}{*}{Age} & \multicolumn{1}{c}{16} & 6\% \\
 & \multicolumn{1}{c}{17} & 6\% \\
 & \multicolumn{1}{c}{18} & 19\% \\
 & \multicolumn{1}{c}{19} & 19\% \\
 & \multicolumn{1}{c}{20} & 19\% \\
 & \multicolumn{1}{c}{21} & 19\% \\
 & \multicolumn{1}{c}{22} & 3\% \\
 & \multicolumn{1}{c}{23} & 3\% \\
 & \multicolumn{1}{c}{24} & 3\% \\
 \midrule
\multirow{6}{*}{Gender} & Cisgender Woman/Girl & 41\% \\
 & Cisgender Man/Boy & 37\% \\
 & Transgender Woman/Girl & 6\% \\
 & Transgender Man/Boy & 6\% \\
 & Genderqueer & 6\% \\
 & Non-binary & 6\% \\
 \midrule
Immigration Status & \multicolumn{1}{c}{Yes} & 15\% \\
\midrule
Human Trafficking Victim Status & \multicolumn{1}{c}{Yes} & 20\%\\
\bottomrule
\end{tabular}
\label{t:youth_demographic}
\end{table}

\begin{table}[]
\centering
\caption{The proportion of RHY requesting certain support services. }
\begin{tabular}{lc}
\toprule
Services & Proportion of RHY \\
\midrule
Mental Health Support & 17\% \\
Medical Support & 50\% \\
Substance Abuse Support & 10\% \\
Crisis 24-Hour Services & 95\% \\
Long-Term Housing Support & 95\% \\
Legal Assistance & 50\% \\
Service Coordination & 95\% \\
Practical Assistance & 50\% \\
Financial Assistance & 50\% \\
Life Skills Support & 50\% \\
Employment Assistance & 50\% \\
Education Assistance & 10\% \\
Childcare Support & 4\%\\
\bottomrule
\end{tabular}
\label{t:youth_needs}
\end{table}

\end{APPENDICES}

%%%%%%%%%%%%%%%%%
\end{document}